\documentclass[a4paper,12pt]{amsart}

\usepackage{amssymb, amscd}
\usepackage{latexsym,epsfig}
\usepackage[all]{xy}
\usepackage{pst-all}
\numberwithin{equation}{section}
\def\today{\ifcase\month\or Jan\or Febr\or  Mar\or  Apr\or May\or Jun\or  Jul\or Aug\or  Sep\or  Oct\or Nov\or  Dec\or\fi \space\number\day, \number\year}

\newcommand{\CC}{\mathbb C}
\newcommand{\EE}{\mathbb E}
\newcommand{\FF}{\mathbb F}

\newcommand{\LL}{\mathbb L}

\newcommand{\PP}{\mathbb P}
\newcommand{\QQ}{\mathbb Q}

\newcommand{\VV}{\mathbb V}

\newcommand{\ZZ}{\mathbb Z}
\newcommand\M[1]{{\mathcal M}_{#1}}

\newcommand\Hy[1]{{\mathcal H}_{#1}}
\newcommand\barM[1]{\overline{\mathcal M}_{#1}}

\newcommand\A[1]{{\mathcal A}_{#1}}

\newcommand\langepijl[1]{\buildrel {#1} \over \longrightarrow}

\numberwithin{equation}{section}

\newtheorem{theorem}{Theorem}[section]
\newtheorem{lemma}[theorem]{Lemma}
\newtheorem{proposition}[theorem]{Proposition}

\newtheorem{conjecture}[theorem]{Conjecture}

\newtheorem{formula}[theorem]{Formula}

\newtheorem{definition-lemma}[theorem]{Definition-Lemma}

\theoremstyle{definition}

\theoremstyle{remark}

\newtheorem{question}[theorem]{Question}

\begin{document}

\title[Curves over Finite Fields and Moduli Spaces]{Curves over Finite Fields and Moduli Spaces}

\author{Gerard van der Geer}
\address{Korteweg-de Vries Instituut, Universiteit van
Amsterdam, Postbus 94248,
1090 GE  Amsterdam, The Netherlands}
\email{g.b.m.vandergeer@uva.nl}


\begin{abstract}
This survey deals with moduli aspects of curves over finite fields. 
We discuss counting points on moduli spaces of curves over finite
fields. Formulas for the number of points on such moduli spaces
lead to modular forms.
In this way counting curves and their points over finite fields
offers a way to obtain information 
on the traces of Hecke operators on modular 
forms. We also discuss stratifications on moduli spaces of curves
and their relevance for curves over finite fields.
\end{abstract}
\maketitle
%%%%%%%%%%%%%%%%%%%%%%%%%%%%%%%%%%%%%%%%
%\centerline{\today}
%%%%%%
\begin{section}{Introduction}
This survey deals with moduli aspects of curves over finite fields. 
A large part of the research on curves in the last century has been on the moduli
of curves. The foundational results on moduli spaces in the 1960s by 
Grothendieck and
Mumford cleared the way for a study of their properties. Before this
it was difficult to study these moduli spaces in positive characteristic.

An attractive aspect of studying curves over a finite field  
is that one can obtain a lot of information 
about them by just counting their 
rational points over the given finite field and its extension fields.
Another nice aspect of studying 
curves over a given finite field  
is that if you fix the genus there are only finitely many of them up to
isomorphism over the given finite field. This suggests that we count these. 
The counting can be done in various ways; for example, one can also
count curves defined over a finite field ${\FF}_q$  of cardinality $q$
together with an $n$-tuple of ${\FF}_q$-rational points. 

Doing the counting for the case
of curves of genus one  leads 
very quickly to the topic of 
modular forms and their Hecke eigenvalues. 
The reason is that counting points of a variety over a finite field 
yields information on the cohomology of that variety. In the case at hand the
variety is a moduli space and modular forms enter naturally in the description
of the cohomology of moduli spaces. Thus we can use the counting of points on curves
over finite fields of cardinality a power of a prime $p$
to study the cohomology of moduli spaces, and to study the traces of
Hecke operators $T_{p^n}$ on spaces of modular forms, and also 
Hecke eigenvalues of modular forms. In this way counting points on curves
over finite fields becomes a heuristic tool in the exploration of 
uncharted terrain.

The interesting thing is that this connection 
can also be used in both directions. 
Knowing traces of Hecke operators provides us with formulas
for the number of rational points on moduli spaces over finite fields.

Moduli of curves admit several stratifications. Examples of these are
the stratification by automorphism group
or by gonality. Apart from these that apply to all characteristics, 
the moduli of
curves in positive characteristic possess stratifications that are special for
given positive characteristic. These stratifications can be quite relevant for
curves over finite fields. There has been a flurry of activity on
these stratifications. We discuss just a few aspects of these 
stratifications and point to several open questions.
\end{section}
\section*{Acknowledgement}
This survey covers mainly developments where I participated or 
contributed. I owe a lot to my collaborators in our joint 
projects through these years and enjoyed enormously the pleasant and
fruitful cooperation.
Thanks are due to Marcel van der Vlugt, Torsten Ekedahl, 
Carel Faber and Jonas Bergstr\"om. 
%%%%%%%%%%%%%%%%%%%%%%%%%%%%%%%%%%%%%%%%
\begin{section}{Moduli Spaces}
To begin at the beginning, assume that we fix non-negative 
integers $g$ and $n$ with $2g-2+n\geq 1$.  A smooth projective geometrically irreducible
curve $C$ of genus $g$
together with $n$ distinct labeled 
points $P_1,\ldots,P_n$ is called an $n$-pointed
curve of genus $g$. It is called stable if the group of automorphisms of $C$
fixing $P_1,\ldots,P_n$ is finite. 
Here a morphism of $(C,P_1,\ldots,P_n)$ to $(D,Q_1,\ldots,Q_n)$
is a morphism $\varphi: C \to D$ with $\varphi(P_i)=Q_i$ for $i=1,\ldots,n$.
The most basic fact here, due to Mumford, 
is the existence for $2g-2+n\geq 1$ 
of a moduli space $\M{g,n}$  of stable $n$-pointed 
curves of genus $g$, \cite{Mumford}.  
It is an irreducible  Deligne-Mumford stack of
dimension $3g-3+n$ defined over $\ZZ$. That we have to deal with stacks and not just
with varieties, is due to the fact that curves can have non-trivial automorphism groups.

The moduli spaces $\M{g,n}$ are in general not complete since curves can degenerate. To
compactify these one introduces the notion of a nodal $n$-pointed curve of genus $g$; it
is a reduced connected proper curve $C$ with finitely many singular points 
which are ordinary double
points such that $C$ has arithmetic genus $g$ and non-singular 
points $P_1,\ldots,P_n$. Again we call such a curve stable
if the automorphism group of $C$ fixing $P_1,\ldots,P_n$ is finite. The basic 
existence result on moduli spaces just mentioned can be strengthened to the fact that there exists
a moduli space $\barM{g,n}$, defined over ${\ZZ}$, 
of stable nodal $n$-pointed curves of genus $g$. It is again a Deligne-Mumford stack
defined over ${\ZZ}$.
The moduli of non-singular stable $n$-pointed curves of genus $g$ form an
open dense part $\M{g,n}$ of $\barM{g,n}$.
This cornerstone theorem is due to Deligne and Mumford (\cite{Deligne-Mumford}) and was established in 1969.

For $g=1$ we need $n\geq 1$ and for $n=1$ we find the moduli space $\M{1,1}$ of elliptic curves.
For $g\geq 2$ we can have $n=0$ and we then just write $\M{g}$ and
$\barM{g}$ for $\M{g,0}$ and $\barM{g,0}$.

Since the moduli space $\M{g,n}$ 
is defined over ${\ZZ}$ we can consider its fibre
$\M{g,n} \otimes {\FF}_p$ in characteristic $p$. These are the moduli spaces
that here we are interested in. We also have their compactifications
$\barM{g,n}\otimes {\FF}_p$.

If $(C,P_1,\ldots,P_n)$ 
is a (smooth) projective curve defined over a finite field 
${\FF}_q$, that is, $P_i\in C({\FF}_q)$ for $i=1,\ldots,n$,
then it defines a ${\FF}_q$-valued point $[C,P_1,\ldots,P_n]$ 
of $\M{g,n}$. But the point
of the moduli space may a priori be defined over a smaller field. In general if
$(C,P_1,\ldots,P_n)$ is defined over a field $L$ that is a Galois extension
of the field $K$ with Galois group $G_{L/K}$, then we can view 
$(C,P_1,\ldots,P_n)$ as a scheme over $K$ and we have an exact sequence
$$
1 \to {\rm Aut}_L((C,P_1,\ldots,P_n)) \to {\rm Aut}_K((C,P_1,\ldots,P_n))
\langepijl{\alpha} G_{L/K} \, .
$$
Then we have:
\begin{lemma}\label{descent} The moduli point 
$[C,P_1,\ldots,P_n]$ can be defined over $K$ if and only if $\alpha$ is 
surjective; the stable $n$-pointed curve $(C,P_1,\ldots,P_n)$ can be defined over $K$ if and only if
$\alpha$ admits a lift.
\end{lemma} 

The second statement is (a variation of) 
a well-known theorem of Weil \cite[Thm 1]{Weil}.
In particular, since in our case of finite fields we may restrict
to the case where $G_{L/K}$ is a cyclic group, 
we see that a ${\FF}_q$-valued point
of $\M{g,n}\otimes {\FF}_p$ can be represented 
by a $n$-pointed curve defined over ${\FF}_q$. 

Thus we can count the number of 
${\FF}_q$-valued points of $\M{g,n}$.
But here the stacky character of the moduli space comes into play.
This stacky aspect means that when we count we have to take into
account the automorphisms of our objects. For example, when dealing
with $\M{1,1}$, the moduli of elliptic curves, we have
$$
\# \M{1,1}({\FF}_q)= \sum_{E} \frac{1}{ \# {\rm Aut}_{{\FF}_q}(E)}
\, ,
$$
where the sum is over all elliptic curves $E=(C,P_1)$ defined over ${\FF}_q$
up to isomorphism over ${\FF}_q$.  

In the case at hand we find $\# \M{1,1}({\FF}_q)=q$. Indeed, there are $q$ 
possibilities for the $j$-invariant of an elliptic curve over ${\FF}_q$.
The justification for this is
that there is another type of moduli space, the coarse moduli space $M_{1,1}$, 
which is a scheme with the property that its $k$-valued points for an 
algebraically closed field $k$ correspond bijectively to the $k$-isomorphism
classes of elliptic curves; this is the $j$-line, where $j$ refers to the
famous $j$-invariant of elliptic curves. For each value of $j \in {\FF}_q$ 
we might have several ${\FF}_q$-isomorphism classes of elliptic curves defined
over ${\FF}_q$, but they contribute in total $1$ in the sum.
In fact, there is the following lemma, 
see \cite[Prop.\ 5.1]{vdG-vdV1992b} for a proof.

\begin{lemma}\label{massgerbe}
Let $C$ be a stable curve defined over ${\FF}_q$. Then we have
$$
\sum_{C'} \frac{1}{\# {\rm Aut}_{{\FF}_q}(C')} = 1\, ,
$$
where the sum is over representatives $C'$ of the ${\FF}_q$-isomorphism classes
contained in the $\overline{\FF}_q$-isomorphism class of $C$.
\end{lemma}  
We thus can ask for $\# \M{g}({\FF}_q)$ for $g\geq 2$. To find this number 
we can make a list of all isomorphism classes 
of curves of the given genus $g$
and determine for each curve in the list 
the order of the automorphism group.

This leads immediately to the question how to find all isomorphism
classes and how to calculate the automorphism groups. Is the answer
always a polynomial in $q$ ? We will deal with these questions in the
next sections. 

Determining the automorphism group of a curve can be difficult.
One remark is that one can avoid calculating the order of the automorphism groups
by considering a family of curves in normal form whose generic member 
has no non-trivial automorphisms such that this family 
is a finite cover of the moduli space.  Then we count in 
this family and divide by the degree of the cover. 

\bigskip

Other moduli spaces that enter here are the moduli spaces of principally
polarized abelian varieties of a given dimension $g$. 
These appear if we study curves via their Jacobians.
Let $\A{g}$ be the moduli space of principally polarized abelian varieties
of dimension~$g$.
This is again an irreducible Deligne-Mumford stack defined over ${\ZZ}$. By associating
to a curve its Jacobian we obtain a map of stacks for $g\geq 2$, the Torelli map,
$$
t=t_g: \M{g} \to \A{g}, \quad [C] \mapsto [{\rm Jac}(C)]
$$
and $t_1: \M{1,1} \langepijl{\sim} \A{1}$. Note that the relative dimension
of $\A{g}$ over ${\ZZ}$ is $g(g+1)/2$ and that of $\M{g}$ is $3g-3$ and
the codimension of the image is $(g-2)(g-3)/2$ for $g\geq 2$.
For $g=2$ the map $t_2$ is an embedding of $\M{2}$ as an open part
of $\A{2}$. But for $g\geq 3$ the map $t_g$ is of stacky degree $2$
onto its image, due to the fact that every abelian variety has an
automorphism of order $2$, while the generic curve has no non-trivial
automorphism.
\end{section}
%%%%%%%%%%%%%%%%%%%%%%%%%%%%%%%%%%%%%%%%
\begin{section}{Counting points of $\M{g,n}$ over finite fields}\label{section-counting}
The first case deals with the moduli spaces $\M{0,n}$ of stable 
$n$-pointed smooth curves of genus $0$. This implies $n\geq 3$.
The coarse moduli space $M_{0,3}$ equals one point and 
the coarse moduli space $M_{0,4}$ equals ${\PP}^1-\{0,1,\infty\}$. 
For $n\geq 3$ we have the formula
$$
\# \M{0,n}({\FF}_q)= \prod_{i=2}^{n-2} (q-i)\, ,
$$
settling the case $g=0$.

For $g=1$ we have $\M{1,1}=\A{1}$ and we find $\# \M{1,1}({\FF}_q)=q$
as observed above. Since we know normal forms for elliptic curves we
can write down a list of all isomorphism classes over ${\FF}_q$. For example,
if the characteristic is not $2$, we can write every elliptic curve 
defined over ${\FF}_q$ as $y^2=f$ with $f\in {\FF}_q[x]$ 
of degree~$3$
with non-vanishing discriminant. 
We find for $q=3$ the following list of isomorphism classes of
elliptic curves over ${\FF}_3$. The first column gives the
polynomial $f$ defining the curve $y^2=f$. The last column, 
that gives the $j$-invariants, illustrates Lemma \ref{massgerbe}.

\smallskip

%\begin{table}
\begin{center}
\begin{tabular}{|| l | c| c|r|}
\hline
$f$ & $\# C({\FF}_3)$ & $1/\# {\rm Aut}_{{\FF}_3}(C)$ & $j$ \\ \hline \hline
$x^3+x^2+1$ & $6$ & $1/2$ & $-1$ \\ \hline
$x^3-x^2-1$ & $2$ & $1/2$ & $-1$ \\ \hline
$x^3+x^2-1$ & $3$ & $1/2$ & $1$  \\ \hline
$x^3-x^2+1$ & $5$ & $1/2$ & $1$  \\ \hline
$x^3+x    $ & $4$ & $1/2$ & $0$  \\ \hline
$x^3-x    $ & $4$ & $1/6$ & $0$  \\ \hline
$x^3-x+1  $ & $7$ & $1/6$ & $0$  \\ \hline
$x^3-x-1  $ & $1$ & $1/6$ & $0$  \\ \hline
\end{tabular}
\end{center}
%\end{table}
We deduce from this table  the frequency list

\begin{center}
\begin{tabular}{|l|r|r|r|r|r|r|r|}
\hline
m & 1 & 2& 3& 4& 5& 6 & 7\\
\hline
freq & 1/6 & 1/2 & 1/2 & 2/3& 1/2 & 1/2 & 1/6\\ \hline
\end{tabular}
\end{center}
where the frequency for given $\# C(k)=m$ is obtained by adding 
the contributions $1/\# {\rm Aut}_k(C)$.

Given this frequency list we know $\# \M{1,n}({\FF}_q)$ for this value of $q$ and
all $n$. Indeed, using the map $\M{1,n} \to \M{1,1}$ we have
$$
\# \M{1,n}({\FF}_q)= \sum_{m} {\rm freq}(m) \binom{m-1}{n-1} (n-1)! \, .
$$
Interpolating the answers for various $q$ 
one finds experimentally for low values
\begin{center}
\begin{tabular}{|l|r|}     \hline
$n$ & $\# \M{1,n}({\FF}_q)$ \\
\hline
$1$ & $q$ \\
$2$ & $q^2$ \\
$3$ & $q^3-1$ \\
$4$ & $q^4-q^2 -3\, q +3$ \\
$5$ & $q^5-5\, q^3-q^2+15\, q-12$\\ 
$6$ & $q^6-15\, q^4+ 25\, q^3 +19\, q^2-80\, q+60 $\\ 
$7$ & $q^7-35\, q^5+125 \, q^4 -126 \, q^3 -155 \, q^2 +490\, q -360$ \\ 
\hline
\end{tabular}
\end{center}
Notice that the degree in $q$ is 
in accordance with the fact that $\dim \M{1,n}=n$.
One may continue and find
$$
\begin{aligned}
\# \M{1,10}({\FF}_q)=
q^{10}-210 \, q^8+2274\, q^7-11655\, q^6+34944\, q^5-62140\, q^4 &\\
+42126\, q^3+89124\, q^2-245664\, q+181440 \, ,  & \\
\end{aligned}
$$
where one may check that $\# \M{1,10}({\FF}_q)$ vanishes for 
$q=2$ and $q=3$.
This answer looks already a bit complicated and raises the question:
\begin{question}\label{vraag1} 
 Is $\# \M{1,n}({\FF}_q)$ always a polynomial in $q$? 
More generally, is $\# \M{g,n}({\FF}_q)$ or
$\# \barM{g,n}({\FF}_q)$ polynomial?
\end{question}
Later we shall see that these experimentally obtained values
given above in the tabel
are correct for all prime powers $q$.

\bigskip

Since Weil and Deligne we know that counting points over finite fields 
of a variety defined over a finite field 
gives information on the cohomology and vice-versa knowledge of the
cohomology tells us about the number of rational points. 

The Lefschetz fixed point theorem connects the number of rational
points on a separated scheme of finite type over $\overline{\FF}_q$,
that is, the number of fixed points of Frobenius, to the trace of 
Frobenius acting on \'etale compactly supported cohomology,
see \cite[Th\'eor\`eme 3.2]{Deligne}. 
By work of Behrend \cite{Behrend}
 this result extends to the setting of Deligne-Mumford stacks.

Question \ref{vraag1} leads to the following question.

\begin{question}\label{vraag2}
What does it mean for the
cohomology that we find polynomials in $q$?
\end{question}

The answer is given by a theorem by
van den Bogaart and Edixhoven \cite{vdB-E2005}.

\begin{theorem} Let ${\mathcal X}$ be a Deligne-Mumford stack
that is proper smooth and of pure dimension $d$ over ${\ZZ}$.
Suppose that for all primes $p$ in a set $S$ of Dirichlet density $1$
there exists a polynomial $P=\sum_{i\geq 0} P_i x^i \in {\QQ}[x]$ such that
$$
\# {\mathcal X}({\FF}_{p^n})=P(p^n) + o(p^{nd/2}) \quad (n \to \infty)\, .
$$
Then $P\in {\ZZ}[x]$ has degree $d$ and satisfies $P(x)=x^d P(1/x)$ and we have
$\#{\mathcal X}({\FF}_{p^n})=P(p^n)$ for all primes $p$ and all $n\geq 1$.
\end{theorem}
The statement says that if there exists a polynomial $P$ 
with rational coefficients 
such that $\lim_{n\to \infty} 
|\#\mathcal{X}({\FF}_{p^n})-P(p^n)| p^{-nd/2}=0$ for enough primes $p$,
then $\# \mathcal{X}({\FF}_{p^n})$ 
is an integral polynomial and of degree $d$.

Another way of expressing it is by saying 
that the cohomology of $\mathcal{X}$ is a polynomial $P({\LL})$
in the Lefschetz motive ${\LL}$. Or, equivalently, 
that $H^i_{\rm et}({\mathcal X}_{\QQ},{\QQ}_{\ell})$ with $\ell\neq p$
vanishes for $i$ odd and equals ${\QQ}_{\ell}(-i/2)^{P_{i/2}}$ for $i$ even. 
Or, one could say  
that geometric Frobenius in ${\rm Gal}(\overline{\FF}_q/{\FF}_q)$ 
acts on $H^i_{\rm et}(\mathcal{X}\otimes \overline{\FF}_p, \overline{\QQ}_{\ell})$ 
for $\ell \neq p$ and even $i$ with eigenvalues $q^{i/2}$. 

\smallskip

But we should bear in mind that, as remarked above, $\M{g,n}$ 
is in general not complete. In order to apply the theorem
just given we need to consider the compactification $\barM{g,n}$.

The complement $\partial \M{g,n}$ of $\M{g,n}$ in $\barM{g,n}$ is
a union of divisors. This complement is
stratified and the strata are quotients by finite groups 
of products of $\M{g',n'}$ for $g'\leq g$
and $n'\leq n+g-g'$. In particular, one sees that even if one is interested in
$\M{g}$ only, the spaces $\M{g,n}$ naturally appear.

\smallskip

There is an action of the symmetric group $\mathfrak{S}_n$ on $\M{g,n}$ and
$\barM{g,n}$. We can then count equivariantly. That is, 
instead of counting fixed points of Frobenius 
$F$ on $\M{g,n}(\overline{\FF}_p)$, we count the
fixed points of $F\circ \sigma$ for $\sigma \in \mathfrak{S}_n$.
This depends only on the cycle type of $\sigma$. 
So we count the number of $\overline{\FF}_q$-isomorphism 
classes of curves together with an $n$-tuple of points
$(P_1,\ldots,P_n)$ on the curve such that it is fixed by 
$F \circ \sigma$.
For example, one can use Lemma \ref{descent} 
to see that the number of fixed points of $F\circ \sigma$ with
$\sigma=(1\, 2)\in \mathfrak{S}_2$ for $\M{g,2}(\overline{\FF}_q)$ 
equals
$$
\sum_{C} \frac{ \# C({\FF}_{q^2})- \# C({\FF}_q)}{\# {\rm Aut}_{\FF_q}(C)}\, ,
$$
where the sum is over all curves $C$ of genus $g$ over ${\FF}_q$ up to
${\FF}_q$-isomorphism.

There is a formula of Getzler-Kapranov (see \cite{G-K}, also \cite{B-T})
that expresses $\# \barM{g,n}({\FF}_q)$
in terms in terms of the $\mathfrak{S}_{n'}$-equivariant counts of
$\M{g',n'}({\FF}_q)$ for $g'\leq g$ and $n'\leq n+g-g'$.

This shows that we can use induction to apply the theorem above 
not only to $\barM{g,n}$, but also to $\M{g,n}$. But this shows
also that one needs the equivariant formulas in order to calculate the
contributions of the boundary strata.

The paper by Getzler and Kapranov deals with the cohomology, but we have 
translated it here in terms of number of rational points.
See also the paper \cite{Diaconu2020} by Diaconu,
who gives (following Getzler-Kapranov) effective formulas
expressing the $\mathfrak{S}_n$-equivariant
Euler characteristics of $\barM{g,n}$ in terms
of the Euler characteristics of $\M{g',n'}$ with the indices
$g',n'$ restricted by 
$\max \{0,3-2g'\} \leq n' \leq 2(g-g')+n$.
%%%%%%%%%%%%%%%%%%%%%%%%%%%%%%%%%%%%%%%%%%%%%%%%%%%%%%%%%%%%%%%%%5555
\end{section}
\begin{section}{Polynomial formulas}
The first equivariant counting for $\M{g,n}$ 
was done by Kisin-Lehrer in \cite{K-L} for $g=0$ and they
applied it to conclude facts about the cohomology of $\M{0,n}$.
They found polynomial functions in $q$ for these equivariant counts
of $\# \M{0,n}({\FF}_q)$.
They showed for example that the alternating representation
of $\mathfrak{S}_n$
does not show up in the cohomology of $\M{0,n}$. 

Getzler did
equivariant counting for $g=1$ in \cite{Getzler1999}. 
In the preceding section we gave a sample of the formulas 
for $\# \M{1,n}({\FF}_q)$  and we add here a sample of those for
$\# \barM{1,n}({\FF}_q)$.

%%%%
\vbox{
\bigskip\centerline{\def\quad{\hskip 0.6em\relax}
\def\quod{\hskip 0.5em\relax }
\vbox{\offinterlineskip
\hrule
\halign{&\vrule#&\strut\quod\hfil#\quad\cr
height2pt&\omit&&\omit&\cr
%\noalign{\hrule}
&$n$ && $\# \barM{1,n}({\FF}_q)$ &\cr
\noalign{\hrule}
&$1$ && $q+1$ &\cr
&$2$ && $q^2+2q+1$ &\cr
&$3$ && $q^3+5q^2+5q+1$ &\cr
&$4$ && $q^4+12q^3+23q^2+12q+1$ &\cr
&$5$ && $q^5+27q^4+102q^3+102q^2+27q+1$ &\cr
&$6$ && $q^6+58q^5+421q^4+756q^3+421q^2+58q+1$ &\cr
} \hrule}
}}

Note the symmetry in the formulas for $\barM{1,n}$ 
displaying Poincar\'e duality.
We will come back to the case $\M{1,n}$ and $\barM{1,n}$ in the next section.

Polynomial formulas for higher genus 
 were first obtained by Getzler for $g=2$ and $n \leq 3$ in \cite{Getzler1998b} 
and then by Bergstr\"om for $g=2$ and $g=3$ 
for some $n$ in \cite{JB1,JB2}. Since these tell us about the cohomology, which is a representation space of $\mathfrak{S}_n$, the answer is
phrased in terms of Schur functions ${\bf s}_{\lambda}$, where
$\lambda$ runs through the partitions of $n$ and these $\lambda$ 
correspond to the irreducible representations of $\mathfrak{S}_n$. 
Recall that the Schur functions form a basis for the symmetric functions, like the elementary or complete symmetric functions. We
refer for the representation theory to \cite{F-H}.
We write the answer as
$$
\sum_{\lambda \vdash n} P_{\lambda}(q) \, {\bf s}_{\lambda}
\quad \text{\rm with} \quad 
P_{\lambda}(q)= \frac{1}{n!} \sum_{\sigma \in \mathfrak{S}_n}
\chi_{\lambda}(\sigma) \, | \M{g,n}^{F\cdot \sigma}| \, ,
$$
with $\chi_{\lambda}$ the character of the irreducible 
irreducible representation determined by $\lambda$ and 
$| \M{g,n}^{F\cdot \sigma}|$ the number of fixed points of 
$F\cdot \sigma$.
To give an idea, here is the equivariant 
answer for $\# \barM{2,4}({\FF}_q)$:
$$
\begin{matrix}
({q}^7+8{q}^6+ 33{q}^5+67{q}^4+67{q}^3+33{q}^2+8{q}+1)\, {\bf s}_4\cr
+(4{q}^6+26{q}^5+60{q}^4+60{q}^3+26{q}^2+4{q})\, {\bf s}_{31} \cr
+(2{q}^6+12{q}^5+28{q}^4+28{q}^3+12{q}^2+2{q})\, {\bf s}_{22} \cr
+(3{q}^5+10{q}^4+10{q}^3+3{q}^2)\, {\bf s}_{211} \cr
\end{matrix}
$$

As one may surmise and we will see, we cannot expect that 
$\# \M{g,n}({\FF}_q)$ is always a polynomial
in $q$. But it happens to be so for small $g$ and $n$. Below we 
summarize a number of cases where
it turns out to be the case. In each of the cases one can normalize an equation
for the curve and then count.

\bigskip

For $g=1$ the $\mathfrak{S}_n$-equivariant version of 
$\# \M{1,n}({\FF}_q)$ 
is polynomial for $1\leq n\leq 10$. 
Similarly, $\mathfrak{S}_n$-equivariant version
of $\barM{1,n}({\FF}_q)$ is polynomial for $1\leq n\leq 10$. 
This is due to Getzler, see \cite{Getzler1999};
we will come back to this in the next section.

For $g=2$  the $\mathfrak{S}_n$-equivariant version of $\# \M{2,n}({\FF}_q)$ 
is polynomial for $0\leq n\leq 9$. The formulas are due to 
Bergstr\"om, \cite{JB2}
for $n\leq 7$, but using the results of Petersen \cite{Petersen} it extends to $n=9$. 
Similarly, this holds also
 $\mathfrak{S}_n$-equivariantly for $\# \barM{2,n}({\FF}_q)$. 

For $g=3$ the $\mathfrak{S}_n$-equivariant version of 
$\# \M{3,n}({\FF}_q)$ 
is polynomial for $0\leq n\leq 7$. The $\mathfrak{S}_n$-equivariant 
version of $\# \barM{3,n}({\FF}_q)$ 
is polynomial for $0\leq n\leq 9$, see \cite{JB1}.

For $g=4$ not much is known. It is expected that $\# \M{4,n}({\FF}_q)$ is polynomial
for $n\leq 3$. The cohomology of $\M{4,0}$ with its Hodge structure was
computed by Tommasi \cite{Tommasi}. These results suggest 
for $\M{4}$ and $\barM{4}$ the expected formulas
$$
\begin{aligned} 
\# \M{4}({\FF}_q) & =q^9+q^8+q^7-q^6\, ,\\
\# \barM{4}({\FF}_q)=q^9+4\, q^8 + 13\, q^7+32\, q^6 & +50\, q^5+ 50\, q^4+32\, q^3+
13\, q^2+4\, q+ 1 \, . \\
\end{aligned}
$$
For $q=2$ this can be confirmed using counts of Xarles \cite{Xarles}.

%%%%%%%%%%%%
\end{section}

%%%%%%%%%%%%%%%%%%%%%%%%%%%%%%%%%%%%%%%%%
\begin{section}{Modular Forms Appear}
As noted above, we cannot expect that $\# \M{g,n}({\FF}_q)$ is always polynomial in $q$
for $g\geq 1$.
This can be illustrated for the case $g=1$, as we will do now.
But instead of considering $\# \M{1,n}({\FF}_q)$, 
to simplify things we look 
at it from a different perspective. 

For each elliptic curve
$E=(C,P_1)$ defined over ${\FF}_q$ 
we have by Hasse $\#E({\FF}_q)= q+1-\alpha-\bar{\alpha}$ with
$\alpha=\alpha(E)$ an algebraic integer with $|\alpha|=\sqrt{q}$. 
Thus we can consider for a non-negative integer $n$
$$
\sigma_n(q):= -\sum_{E} \frac{\alpha^n+\alpha^{n-1}\bar{\alpha}+\cdots + \bar{\alpha}^n }{\# {\rm Aut}_{{\FF}_q}(E)}\, ,
$$
where we sum over all elliptic curves $E$ defined over ${\FF}_q$ up to ${\FF}_q$-isomorphism
and the $\alpha$ depend on $E$.

For example, above we gave 
a frequency list for $p=3$ which immediately provides the values of $\sigma_n(3)$. 
It is easy to see that $\sigma_n(q)=0$ for $n$ odd, since if $\alpha,\bar{\alpha}$ occurs, then $-\alpha,-\bar{\alpha}$ occurs for a twist
with the same factor $1/\#{\rm Aut}_{{\FF}_q}$. 
Working out the frequency lists for
$q=2,3,5,7$ produces the following table.

\begin{center}
\begin{tabular}{|l|r|r|r|r|r|r|r|r|r|r|}
\hline
$n$ & $0$ & $2$ & $4$ & $6$ & $8$ & $10$ & $12$ & $14$ & $16$  \\
\hline \hline
$\sigma_n(2)$ & $-2$ & $1$ & $1$ & $1$ & $1$ & $-23$ & $1$ & $217$ & $-527$ \\ \hline
$\sigma_n(3)$ & $-3$ & $1$ & $1$ & $1$ & $1$ & $253$ & $1$ & $-3347$ & $-4283$ \\ \hline
$\sigma_n(5)$ & $-5$ & $1$ & $1$ & $1$ & $1$ & $4831$ & $1$ & $52111$ & $-1025849$ \\ \hline
$\sigma_n(7)$ & $-7$ & $1$ & $1$ & $1$ & $1$ & $-16743$ & $1$ & $2822457$ & $3225993$ \\ \hline
\end{tabular}
\end{center}
If we subtract the ubiquitous $1$ we recognize that  for $p=2,3,5, 7$
$$
\sigma_{10}(p)= 1+\tau(p)\, ,
$$
where 
$$\Delta=\sum_{n\geq 1} \tau(n) q^n= 
q-24 \, q^2+ 252\, q^3 + \cdots =
q \prod_{m\geq 1}(1-q^m)^{24} \eqno(0)
$$ 
is the well-known normalized cusp form of weight $12$ on 
${\rm SL}(2,{\ZZ})$. And similarly, we recognize $\sigma_{14}(p)$
(resp.\ $\sigma_{16}(p)$) as $1+a(p)$ with $a(p)$  the $p$th 
Fourier coefficient of the 
normalized cusp form $q+\sum_{n\geq 2} a(n)q^n$ that generates
the space $S_{16}$ of cusp forms of weight $16$ (resp.\ $S_{18}$ 
of weight $18$) on ${\rm SL}(2,{\ZZ})$. Such counts were done by Birch \cite{Birch}
in the 1960s.
He gave the formulas equivalent to those of $\sigma_k$ for $2\leq k \leq 10$.

In general we have for even $k>0$
$$
\sigma_k(p)= {\rm Tr}(T_p,S_{k+2}) +1
$$
with $S_k=S_k({\rm SL}(2,{\ZZ}))$ the space of cusp forms of weight $k$ on ${\rm SL}(2,{\ZZ})$
and $T_p$ the Hecke operator at $p$. This is an aspect of the 
Eichler-Shimura-Deligne relation to which we now turn.

The expression $\sigma_k(p)$ calculates cohomological information. The moduli space
$\mathcal{A}_1=\M{1,1}$ carries a local system ${\VV}=R^1\pi_{*}({\QQ}_{\ell})$
of rank $2$ 
with $\pi: \mathcal{X} \to \mathcal{A}_1$ the universal elliptic curve. The fibre of
this local system over $[E]$ is the cohomology $H^1_{\rm et}(E,{\QQ}_{\ell})$.

This local system gives rise for each $k>0$ to a local system
${\VV}_k={\rm Sym}^k({\VV})$. As it turns out, for a prime power $q$ the expression $\sigma_k(q)$
equals the trace of Frobenius $F_q$ on the compactly supported cohomology
of ${\VV}_k$. The contribution of an elliptic curve $E/{\FF}_q$
to $\sigma_k(q)$ is minus the trace of Frobenius $F_q$ on
${\rm Sym}^k(H^1_{\rm et}(E\otimes \overline{\FF}_q,{\QQ}_{\ell}))$.

A classical result of Eichler-Shimura (\cite{Shimura1959}) 
says that for even $k>0$
$$
H^1_c(\mathcal{A}_1({\CC}),{\VV}_k \otimes {\CC})\cong S_{k+2}\oplus \overline{S}_{k+2}
 \oplus {\CC} \, , \eqno(1)
$$
and this displays the mixed Hodge structure on $H^1_c(\mathcal{A}_1({\CC}),{\VV}_k)$.
Deligne showed in 1968 in \cite{Deligne1968}
that this result has an analogue for \'etale
$\ell$-adic cohomology.
Moreover, one can relate Frobenius there to the Hecke operator. The 
result may be summarized as follows.

\begin{theorem}\label{g=1result} Let $k$ be an even positive integer and $p$ a prime. 
Then $\sigma_k(p)$
can be expressed in terms of 
the trace of the Hecke operator $T_p$ on the space of cusp forms $S_{k+2}$ of weight $k+2$ on ${\rm SL}(2,{\ZZ})$ as
$$
\sigma_k(p)= {\rm Tr}(T_p,S_{k+2})+ 1\, .
$$
\end{theorem}
This explains the experimental observation given above.

\smallskip
Thus we know $\sigma_k(p)$ if we know the trace of the Hecke operator $T_p$.
And similarly, for prime powers $q$, 
but the formula is slightly different
as $F_q$ does not correspond exactly to $T_q$. 
But we can turn this around and use counting of elliptic curves over finite 
fields to calculate traces of Hecke operators. That one 
does not encounter this, is because we have a closed formula,
the Eichler-Selberg formula, for the traces of the Hecke operators
on the spaces of cusp forms on ${\rm SL}(2,{\ZZ})$.
%%%%%%%%%%%%
\vskip 1 cm

In order to go back to the (cohomology of the) spaces $\M{1,n}$ that we started with, 
one may use a relation between
$\# \M{1,n}({\FF}_{q})$ and the numbers $\sigma_k(q)$ for $k<n$. 
This relation can be beautifully
expressed by using a formula of Getzler (see \cite[p.\ 200]{Getzler1999})

\begin{proposition} (Getzler's formula)
$$
\begin{aligned}
\frac{\# \M{1,n}({\FF}_q)}{n!} &= \text{\rm residue at $0$ of the formal
expression} \\
& \binom{q-t-q/t}{n} \sum_{k=1}^{\infty} \left( \frac{\sigma_k(q)}{q^{2k+1}t^{2k}-1}\right) (t-\frac{q}{t})\,  dt \, . \\
\end{aligned}
$$
\end{proposition}

We thus see that $\# \M{1,n}({\FF}_{q})$ in general
is not a polynomial in $q$ for $n\geq 11$. And indeed, the modular form
$\Delta$ shows up in the formulas for $n=11$;
$$
\begin{aligned}
\# \M{1,11}({\FF}_p)= & p^{11}-330\, p^9 +4575 \, p^8-30657 \, p^7 +124992 \, p^6 -336820 \, p^5 +584550\, p^4 \\
& -406769\, p^3-865316\, p^2+2437776\, p -1814400 -\tau(p)\, , \\
\# \barM{1,11}({\FF}_p)= &
p^{11}+2037\, p^{10} + 213677 \, p^9 +4577630\, p^8+30215924\, p^7 +
74269967 \, p^6 + \\
& \qquad \qquad 30215924\, p^5 + \cdots + 2037 \, p + 1-\tau(p)\, .\\
\end{aligned}
$$
Note that because of Poincar\'e duality the expression $\# \barM{1,11}({\FF}_p)$
possesses a symmetry.
The complicated formulas for $\# \M{1,n}({\FF}_q)$ also explain 
why we preferred to deal with the function $\sigma_k(q)$ instead of
those for $\# \M{1,n}({\FF}_q)$.

The formulas for $\sigma_k(q)$ and $\# \M{1,n}({\FF}_q)$ are displaying one aspect of
the cohomology of the local systems ${\VV}_k$ on $\mathcal{A}_1$. There is a motivic
form of this that incorporates more aspects. 
Scholl showed in \cite{Scholl} the existence of a motive $S[k+2]$
associated to the space $S_{k+2}$ of cusp forms on 
${\rm SL}(2,{\ZZ})$. Then the Eichler-Shimura-Deligne relation takes the form
for $k>0$
$$
H^1_c(\A{1},{\VV}_k)=S[k+2]+1
$$
and this incorporates both (1) and Deligne's generalization for $\ell$-adic
\'etale cohomology.
 The relation between the Hecke operator $T_p$ and geometric Frobenius $F_p$
then implies 
$$
1+ {\rm Tr}(F_p, S[k+2])=  \sigma_k(p) \, .
$$

%%%%%%%%%%%%%%%%%%%%%%
\end{section}
%%%%%%%%%%%%%%%%%%%%%%%%%%%%%%%%%%%%%%%%
\begin{section}{Genus Two}
For genus $2$ the Torelli map $t_2: \M{2} \to \mathcal{A}_2$ is an embedding
with image the open subset whose complement is $\A{1,1}$, the locus
of principally polarized abelian surfaces that are products of elliptic curves with the product polarization.
It is natural to look for an analogue of the $g=1$ formula
$$
{\rm Tr}(T_p,S_{k+2})=-1+\sigma_k(p) \, .
$$
As we saw in the preceding section, 
the contribution of an elliptic curve $E/{\FF}_q$ to $\sigma_k(q)$ is minus
the trace of Frobenius $F_q$ on 
${\rm Sym}^k(H^1_{\rm et}(E\otimes \overline{\FF}_q,{\QQ}_{\ell}))$.
For an abelian surface $X$ over a finite field the $4$-dimensional vector space 
$H^1_{\rm et}(X,{\QQ}_{\ell})$ (with $\ell$ prime to the characteristic) is provided
with a non-degenerate symplectic pairing
$$
H^1_{\rm et}(X,{\QQ}_{\ell}) \times H^1_{\rm et}(X,{\QQ}_{\ell}) \to
{\QQ}_{\ell}(-1)\, ,
$$
which corresponds to the Weil pairing if one identifies $H^1_{\rm et}(X,{\QQ}_{\ell})$ with the dual of the 
$\ell$-adic Tate module of $X$.
Thus the natural analogue of ${\rm Sym}^k$ is an irreducible representation $R_{a,b}$
of ${\rm Sp}(4,{\QQ})$ with highest weight $a\geq b \geq 0$. 

We recall that the irreducible representations of ${\rm Sp}(4,{\QQ})$
are parametrized by the pairs $(a,b)$ of integers with $a\geq b \geq 0$,
and $R_{a,b}$ occurs in ${\rm Sym}^{a-b}(V) \otimes {\rm Sym}^b(\wedge^2V)$
with $V=R_{1,0}$ 
the standard representation of ${\rm Sp}(4,{\QQ})$. We refer to
\cite{F-H} for the representation theory.

For a smooth projective curve $C$ of genus $2$ over ${\FF}_q$ there exist
by Weil algebraic integers $\alpha_1$ and $\alpha_2$ of absolute value $\sqrt{q}$
such that 
$$
\# C({\FF}_{q^n})=q^n+1-\alpha_1^n-\overline{\alpha}_1^n -\alpha_2^n-\overline{\alpha}_2^n \qquad \text{\rm  for all $n\in {\ZZ}_{\geq 1}$}.
$$
The trace of Frobenius on $R_{a,b}(H^1_{\rm et}(X \otimes \overline{\FF}_q, {\QQ}_{\ell}))$,
with $X$ the Jacobian of $C$, 
can be expressed by certain symmetric expressions (Schur functions) 
in these $\alpha_i$ and $\overline{\alpha}_i$.

Summing this over all principally polarized abelian surfaces $X$ over ${\FF}_q$
up to isomorphism over ${\FF}_q$, including the abelian surfaces that are
products of elliptic curves, yields a number $\sigma_{a,b}(q)$ that 
generalizes the $\sigma_k(q)$
of $g=1$.

Around 2002 Carel Faber and I embarked on a 
program to find the analogue
of Theorem \ref{g=1result} for $g=2$.
What we did (\cite{F-vdG1}) 
was counting curves of genus $2$ over finite fields;
for given field ${\FF}_q$ with $q \leq 37$ (later $q < 200$) we
compiled a frequency list
of possible Weil polynomials
as the curve ran through all curves of genus $2$ up to isomorphism over
 ${{\FF}_q}$
(with factor $1/\# {\rm Aut}_{{\FF}_q}(C)$)
and added the contribution from the degenerate curves (corresponding to products 
of elliptic curves).
Thus we calculated $\sigma_{a,b}(q)$ and then interpolated the outcome by polynomials
in $q$ and known motives, like that of $\Delta$. When this no longer worked
we encountered new modular forms.

The modular forms that are expected to appear here
are Siegel modular forms of degree~$2$.
To explain this notion we 
introduce the Siegel upper half space of degree $g$ by
$$
\mathfrak{H}_g=\{ \tau \in {\rm Mat}(g\times g, {\CC}): \tau=\tau^t,
{\rm Im}(\tau)>0\} \, . 
$$
The symplectic group $\Gamma_g={\rm Sp}(2g,{\ZZ})$ is the automorphism group of the
${\ZZ}$-module of rank $2g$ 
generated by a basis $e_1,\ldots,e_g$ and $f_1,\ldots,f_g$ 
and symplectic form $\langle \, , \, \rangle $ 
with $\langle e_i,e_j\rangle=0=\langle f_i,f_j\rangle$ and 
$\langle e_i,f_j\rangle=\delta_{ij}$, the Kronecker $\delta_{ij}$. 
Using this basis 
an element $\gamma \in {\rm Sp}(2g,{\ZZ})$ can be written as a $2 \times 2$ matrix
of $g\times g$ matrices of integers. An element $\gamma=\left(
\begin{matrix} a & b \\ c & d \\ \end{matrix} \right)$  acts on
$\mathfrak{H}_g$ by $\tau \mapsto (a\tau+b)(c\tau+d)^{-1}$.

We now fix
a finite-dimensional irreducible complex representation 
$$
\rho: {\rm GL}(g) \to W \, .
$$
A Siegel modular form of weight $\rho$ and degree $g>1$ is a holomorphic
map
$f:\mathfrak{H}_g \to W$ with 
$$
f((a\tau+b)(c\tau+d)^{-1})= \rho(c\tau+d) f(\tau) \qquad
 \text{\rm for all} \quad
\left(\begin{matrix} a & b \\ c & d \\ \end{matrix}\right) \in {\rm Sp}(2g,{\ZZ})\, .
$$

For $g=2$ we take as representation 
 the irreducible\ ${\rm GL}(2)$-representation
$$
W={\rm Sym}^j({\rm St}) \otimes \det({\rm St})^{\otimes k}
$$
with ${\rm St}$ the standard representation and $j$ and $k$ integers
with $j\geq 0$. 
We call the weight $(j,k)$.
We have the notion of cusp forms and $S_{j,k}=S_{j,k}(\Gamma_2)$ 
denotes the vector space
of cusp forms of weight $(j,k)$. For $j=0$ we are dealing with 
scalar-valued  Siegel modular forms. We also have an algebra of operators, the Hecke 
algebra with operators $T_p$ for every prime $p$; see for example
\cite{BGHZ} and the references there. Scalar-valued Siegel modular forms of degree
$2$ were studied by Igusa in the 1960s \cite{Igusa1962,Igusa1967}. 
Igusa determined the ring of scalar-valued modular forms of degree $2$.
Satoh was one of the first
to study vector-valued ones \cite{Satoh}.

Let us note that 
one may also express Siegel modular forms of degree $2$ 
as sections of a bundle.
If ${\EE}$ is the Hodge bundle on $\A{2}$ whose fibre over $[X]$
is the $2$-dimensional space $H^0(X,\Omega^1_X)$, 
then a Siegel modular
form of weight $(j,k)$ can be viewed as a section of 
${\rm Sym}^j({\EE}) \otimes \det({\EE})^{\otimes k}$.

\smallskip

In 2002 using the frequency counts of Weil polynomials of genus $2$ curves
over finite fields
Carel Faber and I found experimentally for $(a,b)\neq (0,0)$ and $a+b$ even,
a formula 
that expresses the trace of the Hecke operator $T_p$  on the space $S_{a-b,b+3}$ 
of cusp forms of weight $(a-b,b+3)$
in terms of the expressions $\sigma_{a,b}(p)$ obtained by counting. 
The formula has the following form, see \cite{F-vdG1}.

\begin{formula} \label{g=2result} 
Let $a,b$ be non-negative integers with $a+b$  even and positive. 
The trace of the Hecke operator $T_p$ on the 
space $S_{a-b,b+3}$ of degree $2$ Siegel modular cusp forms of weight 
$(a-b,b+3)$ can be expressed in the counting function $\sigma_{a,b}(p)$
by
$$
{\rm Tr}(T_p,S_{a-b,b+3})=\sigma_{a,b}(p)+c_{a,b}(p),
$$
where the `correction term' $c_{a,b}(p)$ equals
$$
s_{a-b+2}-s_{a+b+4} \, \sigma_{a-b+2}(p) \, p^{b+1} + \begin{cases} 
\sigma_{b+2}(p) \cr 1-\sigma_{a+3}(p)\end{cases}
$$
with $s_k=\dim S_k({\rm SL}(2,{\ZZ}))$. 
\end{formula}
This term $c_{a,b}(p)$ (and
more generally a term $c_{a,b}(q)$) is
the analogue of the term $1$ for $g=1$. Note that the correction
involves only $g=1$ data.

This conjectured formula \ref{g=2result} for $c_{a,b}(q)$ 
was later proved by Weissauer (2009) \cite{Weissauer}.
For this one needs the Eisenstein cohomology; we refer to
\cite{vdG2011} and to
Harder \cite{Harder} and Petersen \cite{Petersen} for the
local systems of irregular highest weight.

Counting points over finite fields thus 
provides a very efficient way to calculate traces 
of Hecke operators on cusp forms of degree $2$.
If you have counted for one prime $p$ you know ${\rm Tr}(T_p,S_{j,k})$
{\bf for all $j,k$} with $k\geq 3$. 
Similarly for prime powers $q$.

We illustrate this with two examples where $\dim S_{j,k}=1$; 
then the trace of $T_p$ equals the eigenvalue. 
Eigenvalues for $T_p$
are denoted by $\lambda(p)$. 
One can calculate eigenvalues of a modular form which is an
eigenform under the Hecke algebra if one knows the Fourier
expansion. In general it is difficult to give such Fourier
expansions and also very laborious to calculate the eigenvalues.
Except for very few cases, such eigenvalues were not known 
for vector-valued forms or for scalar-valued forms of larger weight.

The two examples are illustrated by tables of Hecke eigenvalues.
The first table gives eigenvalues $\lambda(p)$ of $T_p$ on  a generator of 
$S_{8,8}$.

\vbox{
\bigskip\centerline{\def\quad{\hskip 0.6em\relax}
\def\quod{\hskip 0.5em\relax }
\vbox{\offinterlineskip
\hrule
\halign{&\vrule#&\strut\quod\hfil#\quad\cr
height2pt&\omit&&\omit&\cr
&$p$&& $\lambda(p)$  on $S_{8,8}$&\cr
height2pt&\omit&&\omit&\cr
\noalign{\hrule}
height2pt&\omit&&\omit&\cr
&$2$&&$2^6\cdot 3 \cdot 7$&\cr
&$3$&&$-2^3\cdot 3^2 \cdot 89$&\cr
&$5$&&$-2^2\cdot 3 \cdot 5^2 \cdot 13^2 \cdot 607$&\cr
&$7$&&$2^4\cdot 7 \cdot 109 \cdot 36973$& \cr
&$11$&&$2^3\cdot 3 \cdot 4759 \cdot 114089$& \cr
&$13$&&$-2^2\cdot 13 \cdot 17 \cdot 109 \cdot 3404113$& \cr
&$17$&&$2^2\cdot 3 \cdot 17 \cdot 41 \cdot 1307 \cdot 168331$& \cr
&$19$&&$-2^3\cdot 5 \cdot 74707 \cdot 9443867$& \cr
height2pt&\omit&&\omit&\cr
} \hrule}
}}

The next table deals with the first non-zero scalar modular form 
of odd weight $35$. This form $\chi_{35} \in S_{0,35}$ is one of
the generators of the ring of scalar-valued Siegel modular forms
of degree $2$ as described by Igusa \cite{Igusa1962}.

%\begin{footnotesize}
\smallskip
\vbox{
\bigskip\centerline{\def\quad{\hskip 0.6em\relax}
\def\quod{\hskip 0.5em\relax }
\vbox{\offinterlineskip
\hrule
\halign{&\vrule#&\strut\quod\hfil#\quad\cr
height2pt&\omit&&\omit&\cr
&$p$&& $\lambda(p)$  on $S_{0,35}$&\cr
height2pt&\omit&&\omit&\cr
\noalign{\hrule}
height2pt&\omit&&\omit&\cr
&$2$&&$-25073418240$&\cr
&$3$&&$-11824551571578840$&\cr
&$5$&&$9470081642319930937500$&\cr
&$7$&&$-10370198954152041951342796400$&\cr
&$11$&&$-8015071689632034858364818146947656$&\cr
&$13$&&$-20232136256107650938383898249808243380$&\cr
&$17$&&$118646313906984767985086867381297558266980$&\cr
&$19$&&$2995917272706383250746754589685425572441160$&\cr
&$23$&&$-1911372622140780013372223127008015060349898320$&\cr
&$29$&&$-2129327273873011547769345916418120573221438085460$&\cr
&$31$&&$-157348598498218445521620827876569519644874180822976$&\cr
&$37$&&$-47788585641545948035267859493926208327050656971703460$&\cr
height2pt&\omit&&\omit&\cr
} \hrule}
}}
%\end{footnotesize}
Note that the parabolic shape of the figures in this diagram 
nicely reflects Deligne's
result on the absolute values of the eigenvalues of Frobenius
\cite{Deligne1968,Deligne1974}.

For $g=1$ knowing the Hecke eigenvalues of a normalized eigenform means knowing the
Fourier expansion of the form. That is no longer the case for $g>1$.

These results on the Hecke eigenvalues of genus $2$ forms 
stimulated Harder to formulate a precise conjecture 
on congruences between elliptic
modular forms and Siegel modular forms, see \cite{BGHZ}. 
For a long time he had suspected that there should be such conguences.
Our work provided strong evidence for
this conjecture.
For me such results beautifully show the use of counting points 
over finite fields.
\bigskip

With Bergstr\"om and Faber we extended this in \cite{BFvdG} 
to $g=2$ and level $2$
(and $p\neq 2$) by taking into account ramification points of the genus $2$
curve $y^2=f$ for $f$ of degree $6$. 
The symmetric group $\mathfrak{S}_6$ acts and we can count equivariantly. 
The equivariant
formulas are somewhat complicated. These were later proved to be correct
by R\"osner \cite{Roesner} using the theory of automorphic representations.

Results on traces of Hecke operators obtained by counting points of
curves over finite fields
can be found on the website {\tt smf.compositio.nl}.

\smallskip

As a final remark, we point out that as soon as cusp forms appear we will not have polynomial formulas for $\# \M{2,n}({\FF}_q)$. In particular, for $n=10$ 
and $\mathfrak{S}_n$-representation 
$\lambda=[1^{10}]$ we see modular forms appear. The dimension of
$S_{j,k}$ grows fast, for fixed $j$ cubically in $k$.

\end{section}
%%%%%%%%%%%%%%%%%%%%%%%%%%%%%%%%%%%%%%%%
\begin{section}{Genus Three}
When one goes from genus $2$ to $3$ the complexity increases.
One aspect of this is 
that for genus $3$ the Torelli map 
$t: \M{3} \to \A{3}$ is a morphism of stacks of degree $2$.
Indeed, every principally polarized abelian variety $X$ of dimension $3$ has an automorphism
$-1_X$ of order $2$ which acts by sending an element to its inverse in the group, 
but the generic curve of genus $3$ does not have such an automorphism. 
A hyperelliptic curve has such an involution and it induces the automorphism $-1_X$ on its
Jacobian. We thus can interpret $\M{3} \to \A{3}$ 
as a double cover ramified along the (closure of)
the hyperelliptic locus. 

Igusa showed in \cite{Igusa1967} that the closure of the locus of hyperelliptic Jacobians in
$\A{3}$ is the zero divisor of a Siegel modular cusp form $\chi_{18}$ of degree $3$ 
and weight $18$. One can view $\chi_{18}$ as a section of the line bundle
${\det}(\EE)^{\otimes 18}$ with ${\EE}$ the rank $3$ 
Hodge bundle on $\barM{3}$.
It allows us to view $\M{3}$ as a double cover of $\A{3}$ obtained by taking a square
root $\chi_9$ of $\chi_{18}$. Ichikawa showed in 1995 that $\chi_9$ is a Teichm\"uller
modular form of weight $9$, 
that is, a section of ${\det}({\EE})^{\otimes 9}$ 
on $\barM{3}$, \cite{Ichikawa}. An algebraic way to construct it
is by observing that there is a natural morphism of locally free 
sheaves of rank $6$  
$$
{\rm Sym}^2({\EE}) \langepijl{} \pi_*(\omega_{\mathcal{C}/\M{3}}^{\otimes 2})
$$
with $\pi: \mathcal{C} \to \M{3}$ the universal curve, obtained
by multiplying differential forms. This morphism
is an isomorphism outside the hyperelliptic locus; taking the determinant
gives a map $\det({\EE})^{\otimes 4} \to \det({\EE})^{\otimes 13}$, hence a section of
$\det({\EE})^{\otimes 9}$ that vanishes on the hyperelliptic locus.

When suitably normalized this form gives the obstruction of an
indecomposable principally polarized abelian threefold $X$ for being a Jacobian of a curve;
it is expressed by saying that
$\chi_{18}$ assumes a square value at $[X]$. We refer to the paper by Ritzenthaler 
\cite{Ritzenthaler} in the Serre book \cite{Serre}. 

Since there are no scalar-valued Siegel modular forms (on the full group ${\rm Sp}(6,{\ZZ})$) of
odd weight, the existence of the Teichm\"uller form $\chi_9$ shows 
that $\M{3}$ carries more cohomology than $\A{3}$.

In joint work with Bergstr\"om and Faber \cite{BFvdG2} 
we found the generalization to $g=3$ 
of the formulas for $g=1$ and $g=2$.
For $X$ a principally polarized abelian threefold over a finite field ${\FF}_q$ 
we consider $H^1(X,{\QQ}_{\ell})$ with $\ell$ prime to $q$.
For varying $X$ it defines a local system ${\VV}$ of rank $6$ on $\A{3}$. The space
$H^1(X,{\QQ}_{\ell})$ carries a non-degenerate symplectic pairing
$$
H^1(X,{\QQ}_{\ell})\times H^1(X,{\QQ}_{\ell}) \to {\QQ}_{\ell}(-1) 
\, .
$$
For each irreducible representation  $R_{a,b,c}$ of ${\rm Sp}(6,{\QQ})$ 
of highest weight $(a,b,c)$ with $a\geq b \geq c$ 
we can form a local system ${\VV}_{a,b,c}$ 
with fibre $R_{a,b,c}(H^1_{\rm et}(X,{\QQ}_{\ell}))$.

The trace of Frobenius on the cohomology 
$$
R_{a,b,c}(H^1_{\rm et}(X\otimes \overline{\FF}_q,{\QQ}_{\ell}))
$$
for $\ell \neq p$ is a symmetric expression (Schur polynomial) 
in the roots of the Weil polynomial of $X$.
Summing this over all $X$ up to isomorphism over ${{\FF}_q}$ gives a number 
$\sigma_{a,b,c}(q)$. This number generalizes $\sigma_a(q)$ 
for $g=1$ and $\sigma_{a,b}(q)$ for $g=2$.

For $g=1$ and $g=2$ we saw the formulas 
$$
{\rm Tr}(T_p,S_{a+2}(\Gamma_1))=\sigma_{a}(p)-1, \quad
{\rm Tr}(T_p,S_{a,b}(\Gamma_2))=\sigma_{a,b}(p)+c_{a,b}(p)\, ,
$$
with the sigmas obtained by counting 
and $-1$ and $c_{a,b}(p)$ a correction term.
Now we deal with Siegel modular forms of degree $3$.
We found experimentally a formula for the trace of the Hecke
operator $T_p$ on the space of cusp forms $S_{i,j,k}=S_{i,j,k}(\Gamma_3)$.
Note that the weight of a modular form on 
$\Gamma_3={\rm Sp}(6,{\ZZ})$ is now given by a triple
$(i,j,k)$. Scalar-valued modular forms correspond to $i=j=0$.

The formula of \cite{BFvdG2} takes the form
$$
{\rm Tr}(T_p, S_{a-b,b-c,c+4})= \sigma_{a,b,c}(p) + c_{a,b,c}(p)\, ,
$$
where $\sigma_{a,b,c}(p)$ is the trace of Frobenius on the compactly supported
cohomology of the local system ${\VV}_{a,b,c}$ over $\A{3}\otimes{\FF}_p$
and the correction term $c_{a,b,c}(p)$ is given by
$$
\begin{matrix}
-\sigma_{a+1,b+1}(p)+\sigma_{a+1,c}(p)-\sigma_{b,c}(p) \\
-c_{a+1,b+1}(p)\sigma_{c+2}(p)+ c_{a+1,c}(p)\sigma_{b+3}(p) -c_{b,c}(p)\sigma_{a+4}(p) \\
+c_{a+1,b+1}(p)-c_{a+1,c}(p)+c_{b,c}(p) \, . \\
\end{matrix}
$$
The correction term is expressed in genus $2$ and $1$ terms and
its form was suggested by formulas in \cite{vdG2011}.

So far, the formula has not yet been proved. But it works perfectly.
That is, there is overwhelming evidence 
that it produces indeed the right values of 
the traces of the Hecke operators on the spaces of cusp forms. 
Not only that, it also allows us to see the contributions of `lifts',
modular forms of degree $1$ and $2$. After identifying these, 
one is left
with the genuine Siegel modular forms of degree $3$ that correspond to
$8$-dimensional Galois representations in the cohomology of a local
system on $\A{3}$. Genuine means that these modular forms do not
belong to $g=1$ and $g=2$. Details can be found in \cite{BFvdG2}.

We refer to the website {\tt smf.compositio.nl} where 
results on the traces of Siegel modular forms of degree $3$ thus obtained 
can be found.
So counting points on curves over finite fields tells us about Siegel modular forms. 
Here are two examples. The spaces $S_{6,3,6}$ and $S_{2,4,8}$ are $1$-dimensional
and the next tables give the eigenvalues $\lambda(p)$ of the Hecke
operators $T_p$ on these spaces.

\smallskip
\vbox{
\bigskip\centerline{\def\quad{\hskip 0.6em\relax}
\def\quod{\hskip 0.5em\relax }
\vbox{\offinterlineskip
\hrule
\halign{&\vrule#&\strut\quod\hfil#\quad\cr
height2pt&\omit&&\omit&&\omit&&\omit \cr
&$p$&& $\lambda(p)$  on $S_{6,3,6}$&&$\lambda(p)$ on  $S_{4,2,8}$& \cr
height2pt&\omit&&\omit&&\omit&&\omit \cr
\noalign{\hrule}
height2pt&\omit&&\omit&&\omit&\cr
&$2$ && $0$&&$9504$& \cr
&$3$ && $-453600$&&$970272$& \cr
&$4$ && $10649600$&&$89719808$& \cr 
&$5$ && $-119410200$&&$-106051896$& \cr
&$7$ && $12572892800$&&$112911962240$& \cr
&$8$ && $0$&&$1156260593664$& \cr
&$9$ && $-29108532600$&&$5756589166536$& \cr
&$11$ &&$ -57063064032$&&$44411629220640$& \cr
&$13$ && $-25198577349400$&&$209295820896008$& \cr
&$16$ && $341411782197248$&&$-369164249202688$& \cr
&$17$ && $-107529004510200$&&$1230942201878664$& \cr
&$19$ && $1091588958605600$&&$51084504993278240$& \cr
height2pt&\omit&&\omit&&\omit&\cr
} \hrule}
}}

If we go to $\M{3}$ instead there is more cohomology. It is known that the cohomology
of local systems on $\A{3}$ can be described in terms of Siegel modular forms (of degree
$\leq 3$). But for $\M{3}$ it is not known what automorphic forms show up in the
cohomology of local systems on $\M{3}$. 
It is a mystery which modular
forms or motives we will encounter. But we can use counting points over finite fields to try to
explore the first cases. For the cohomology of $\M{3}$ we know that we are dealing with
motives or modular forms `of level one', or in other words, with
everywhere good reduction. 

An example is provided by the local system ${\VV}_{11,3,3}$ whose fibre for a curve $C$
is the irreducible representation of highest weight $(11,3,3)$ on 
$H^1_{\rm et}(C,{\QQ}_{\ell})$, see \cite[p.\ 34]{CFvdG}. 
We see a $6$-dimensional motive of weight $23$
appearing in the cohomology 
with Hodge degrees $0,5,9,14,18,23$. 
It cannot come
from a Siegel modular form of degree $3$. It should correspond to a
Teichm\"uller modular form of weight $(8,0,7)$, that is, a section of 
${\rm Sym}^8({\EE}) \otimes {\det}({\EE})^{\otimes 7}$ on $\barM{3}$.
We can construct such a Teichm\"uller modular form. It is fascinating to be able to 
explore this cohomology and these motives simply by counting points on curves
over finite fields. 

\end{section}
%%%%%%%%%%%%%%%%%%%%%%%%%%%%%%%%%%%%%%%%
\begin{section}{Other Cases}
The fact that a moduli space of curves like $\M{g}$ for $g=1,2,3$
maps generically finitely to a Shimura variety helps a
lot for understanding the cohomology of local systems, and hence understanding
the behavior of quantities 
like $\# \M{g,n}({\FF}_q)$. Things become much harder 
if one does not know a priori which
modular forms or motives may show up.  
The case of $\M{4}$ and $\M{4,n}$ illustrates this. We do not 
know what motives or modular forms we can expect. 
Again, for the cohomology of $\M{g,n}$ we know that we are dealing with
motives or modular forms `of level one', or in other words,
with everywhere good reduction.
As far as I know, nobody has an idea of the nature of
 the zeta function of $\M{g}$ for general $g$. 
Results of 
Chenevier, Lannes and Renard 
on motives of level $1$ and small rank
\cite{Chenevier-Lannes,Chenevier-Renard} can help
in identifying the results for low values of $g$ and $n$.

Given the absence of knowledge what answer to expect in general,
it is natural to first
look at cases where moduli spaces of curves are closely related to arithmetic
quotients. In 1964 Shimura gave a list in \cite{Shimura1964}
of arithmetic quotients of the ball
that appear as moduli spaces of curves. Rohde extended this list 
in \cite{Rohde} and Moonen proved his list is complete \cite{Moonen}. 
In these cases one can count and hope to be able to interpret
the result in terms of modular forms for the arithmetic group in question. 

One case of Shimura's list deals with the family of Picard curves.
These are curves of genus $3$ that are cyclic Galois covers of degree $3$ of
the projective line; these can be given by an equation
$y^3=f$
with $f$ a polynomial of degree $4$ with non-zero discriminant. Picard 
studied these curves in the late 19th century. The Jacobians of these
curves have multiplication by the ring of integers $O_F$ 
of $F={\QQ}(\sqrt{-3})$
and the moduli space is a so-called Picard modular surface.
In this case this modular surface 
is a ball quotient associated to a discrete subgroup
of the group of unitary similitudes ${\rm U}(2,1,{\QQ}(\sqrt{-3}))$.

In a long term project with Bergstr\"om (\cite{B-vdG}) 
we have counted points on such curves over
finite fields and using this we arrived heuristically  
at a complete formula for the traces of Hecke operators 
on the corresponding modular forms in terms of these counts
over finite fields. 

The formula expresses the trace of a Hecke operator $T_{\nu}$
for $\nu$ a prime of $O_F$ of norm congruent to $1$ modulo $3$
on a space of cusp forms in terms of counts on curves over finite fields;
more precisely, it expresses it in terms of  an analogue of the functions
$\sigma_k$ and $\sigma_{a,b}$ 
obtained by counting that we saw above, 
and a correction term.

But the results provided further insights.
A careful analysis of the results allowed us to distill from this
the traces of the Hecke operators on the spaces of genuine Picard modular
forms, that is, the modular forms that correspond to Galois representations
of degree $3$.  

The results of \cite{B-vdG} 
are conjectural, but as in the case $g=3$,
it works perfectly. Here the weight of a modular form is
a pair $(j,k)$ of non-negative integers, the case $j=0$
corresponding to scalar-valued modular forms. 

We counted for primes $p\equiv 1 (\bmod \, 3)$ 
with $p\leq 43$. 
For these primes we can predict the trace of the Hecke
operator $T_{\nu}$, for $\nu$ a prime of $O_F$ with
norm $p$, on the space of cusp forms of 
given weight $(j,k)$ and also on the space of 
genuine cusp forms of given weight $(j,k)$ {\bf for any weight}
$(j,k)$ with $k\geq 3$.

We illustrate this by one example, but must refer to 
\cite{B-vdG} for details. It concerns a $1$-dimensional 
space of cusp forms of weight $(0,33)$ with a 
character of order $2$.
The Hecke eigenvalues of a generator of this space, 
that we found, agree with a 
congruence modulo  $17093$ of the type of congruences 
predicted by Harder (see Harder's contribution to \cite{BGHZ}). 

If $p$ is a prime with $p\equiv 1 \bmod 3$ then $p$
splits in $O_F$ as 
$p=\nu_p \bar{\nu}_p$
with $\nu_p\equiv \bar{\nu}_p 
\equiv 1 (\bmod 3)$.  
The Hecke eigenvalues $\lambda_{\nu_p}$ that we found are given 
in the next table. There $\rho$ denotes a third root of $1$.

\begin{footnotesize}
\smallskip
\vbox{
\bigskip\centerline{\def\quad{\hskip 0.6em\relax}
\def\quod{\hskip 0.5em\relax }
\vbox{\offinterlineskip
\hrule
\halign{&\vrule#&\strut\quod\hfil#\quad\cr
height2pt&\omit&&\omit&\cr
&$p$ && $\lambda_{\nu_p}$ & \cr
\noalign{\hrule}
& $7$   && $-50515470688848 + 19722585570921 \rho$ & \cr
& $13$  && $-641186317588670376 - 28497381958498509 \rho$ & \cr
& $19$  && $-207202261228535219325 - 223900464575892946149 \rho$ & \cr
& $31$  && $-72536002932093668516175 - 511708107362090202586656 \rho$ & \cr
& $37$  && $ -15066567237821284922757576 - 12800018433999723562677897 \rho$ & \cr
& $43$  && $2263923296934966075769869 - 61311985796827137336770952 \rho $ & \cr
} \hrule}
}}
\end{footnotesize}
These eigenvalues $\lambda_{\nu_p}$ satisfy the congruence
$$
\lambda_{\nu_p}\equiv \bar{\nu_p}^{32}+(p^{31}+1)\nu_p  \, 
(\bmod \,  17093) \, ,
$$
as the reader may verify.
We refer to \cite[Section 14]{B-vdG}. 

\vskip 1 cm

Besides this case we intend to treat other cases
of Shimura's list.
But one may decide to go further into unknown territory 
and count for other families
of curves not closely related to Shimura varieties. 
One obvious case are the families $\Hy{g}$ of 
hyperelliptic curves of given $g$ and their variants $\Hy{g,n}$. Here
the case of characteristic $2$ plays a special role. Bergstr\"om gives
in \cite{JB1} recursive formulas in the genus 
for the equivariant number of points of $\Hy{g,n}$ over a fixed finite field.
This reduces the problem for fixed $n$ to the cases of low genera.
We refer to \cite{JB1} for more details.
\end{section}
%%%%%%%%%%%%%%%%%%%%%%%%%%%%%%%%%%%%%%%%
\begin{section}{Stratifications}
The moduli spaces $\M{g}$ and $\A{g}$ admit stratifications. Some of these
stratifications work in all characteristics, other ones are special to
positive characteristic. 

Maybe the best known stratification is by automorphism group. Here we
see considerable differences between characteristic zero and 
positive characteristic.
The well-known Hurwitz bound on the order of the automorphism group 
of a curve in
characteristic $0$ is no longer true in positive characteristic. 

For some low genera we know explicit stratifications of 
$\M{g}({\CC})$ by automorphism group.
It would be nice to have such explicit stratifications also 
for $\M{g}({\overline{\FF}}_p)$ and for $\M{g}({\FF}_q)$. 

\smallskip 

Another example of a stratification of $\M{g}$ is given by gonality. 
Here gonality of a curve $C$ over a field $k$ 
is the smallest degree of a morphism $C\to {\PP}^1$ over $k$. 

The gonality of a curve $C$ defined over ${\FF}_q$ 
 puts clear restrictions on the number of ${\FF}_q$-rational points
of $C$. 
Thus gonality may be relevant for the study of the invariant $N_q(g)$,
the maximum number
of rational points on a smooth projective 
curve of genus $g$ over ${\FF}_q$, that
plays such an important role in Serre's lectures notes \cite{Serre}.

In \cite{vdG-ECM} I asked the question: {\sl what is the maximum number of
rational points on a curve of genus $g$ and gonality $\gamma$
defined over ${\FF}_q$?} It suggests to study a variant of $N_q(g)$; 
namely the invariant
$N_q(g,\gamma)$: 
the maximum number of ${\FF}_q$-rational points on a (smooth
projective) curve over ${\FF}_q$ of genus $g$ and gonality $\gamma$.

Recently, this question was taken up 
by X.\ Faber and Grantham for small $q$ and $g$ in
\cite{F-G1,F-G2}. 
It is well-known that the gonality over a finite field is at most $g+1$,
a result due to F.K.\ Schmidt. 
Moreover, if $\# C({\FF}_q)>0$ then $\gamma(C) \leq g$. 
Indeed, for $g\geq 2$ if $P\in C({\FF}_q)$
then $h^0(K-(g-2)P) \geq 2$ by Riemann-Roch, providing a linear system of
degree $g$.

But if $\#C({\FF}_q)=\emptyset$, then the gonality can be $g+1$. For example,
the curve $C$ of genus $3$ defined over ${\FF}_2$ given by
$$
x^4+y^4+z^4+x^2y^2+x^2z^2+y^2z^2+x^2yz+xy^2z+xyz^2=0
$$
in ${\PP}^2$ has no rational points and it is not difficult to show that it has
gonality $4$.
Here is a table with some of their results (taken from \cite{F-G2}).
The cross $\times$ indicates the absence of curves over ${\FF}_q$
with given $(g,\gamma)$.

\begin{center}
\begin{tabular}{|l|r|r|r|r|r|}
\hline
$g$ & $\gamma$ & $N_2(g,\gamma)$ & $N_3(g,\gamma)$ & $N_4(g,\gamma)$  \\
\hline
$3$ & $2$ & $6$ & $8$ & $10$ \\
$3$ & $3$ & $7$ & $10$ & $14$ \\
$3$ & $4$ & $0$ & $0$ & $0$ \\
\hline
$4$ & $2$ & $6$ & $8$ & $10$ \\
$4$ & $3$ & $8$ & $12$ & $15$ \\
$4$ & $4$ & $5$ & $10$ & $13$ \\
$4$ & $5$ & $0$ & $0$ &  {\color{red}$\times$} \\
\hline
$5$ & $2$ & $6$ & $8$ & $10$ \\
$5$ & $3$ & $8$ & $12$ & $15$ \\
$5$ & $4$ & $9$ & $13$ & $17$ \\
$5$ & $5$ & $3$ & $4$ &  $5$ \\
$5$ & $6$ & {\color{red}$\times$} & {\color{red}$\times$} & {\color{red}$\times$} \\
\hline
\end{tabular}
\end{center}

Faber-Grantham conjecture in \cite{F-G2} 
that an optimal curve has gonality at most $\lfloor \frac{g+3}{2} \rfloor$ and that for fixed $\gamma \geq 2$ and fixed $q$ 
and for $g$ large one has 
$N_q(g,\gamma)=\gamma(q+1)$. Together with Howe they show for $g\geq 5$ in
\cite{F-G-H} that over a finite field gonality is at most $g$.

In a recent paper \cite{Vermeulen} Floris Vermeulen shows by a construction
of curves in a toric variety the following result. 

\begin{theorem} For fixed $q$ and $\gamma$ with 
$\gamma\leq q+1$  we have  $\lim_{g \to \infty} N_{q}(g,\gamma)=\gamma \, (q+1)$.
\end{theorem}

\end{section}
%%%%%%%%%%%%%%%%%%%%%%%%%%%%%%%%%%%%%%%%
\begin{section}{Characteristic $p$ stratifications}
There are stratifications on the moduli of curves and abelian varieties
that are special to positive characteristic. The first well-known case where 
this phenomenon appears is the case of elliptic curves. 
Elliptic curves in characteristic $p>0$ come into two sorts: 
ordinary or supersingular.
A formula of Deuring gives the number of supersingular curves.
The number of isomorphism classes of supersingular elliptic
curves over $\overline{\FF}_p$ equals
$$
h_p= \frac{p-1}{12}+ \left( 1-({{-3}\over {p}})\right)\frac{1}{3}+
\left( 1- ({{-4}\over {p}})\right) \frac{1}{4}\, .
$$
But a stacky interpretation gives the more elegant formula
$$
\sum \frac{1}{\# {\rm Aut}(E)} = \frac{p-1}{24} \, ,
$$
where the summation is over all isomorphism classes of
supersingular elliptic curves defined over
$\overline{\FF}_p$.
This stratification on $\A{1}\otimes {\FF}_p$ generalizes for higher $g$ 
to a stratification on $\A{g}\otimes {\FF}_p$ in two ways.

\begin{enumerate}
\item{} The Newton Polygon stratification (NP). 

\item{} The Ekedahl-Oort stratification (E--O). 

\end{enumerate}

The simplest examples of strata in both stratifications
 are the $p$-rank strata.
If $X/k$ is an abelian variety over a field $k$
of characteristic $p>0$ and $X[p]$ denotes the kernel of multiplication
by $p$ on $X$ with $\# X[p](\overline{k})=p^f$ then
the $p$-rank of $X$ is $f$. The (closed) $p$-rank strata are
$$
V_f=\{ [X] \in \A{g}(\overline{k}): f(X) \leq f \}\, .
$$
Koblitz showed already in 1975 in \cite[Thm.\ 7]{Koblitz} that this
gives a stratification with ${\rm codim}(V_f)=g-f$
in $\A{g}\otimes {\FF}_p$.

The study of these stratifications was pursued by Oort
and has become a very active area of research. I summarize
a few salient features before discussing their relevance for
curves over finite fields.

The first one, the Newton Polygon stratification,
extends the $p$-rank stratification and
 was introduced by Grothendieck and Katz.  They showed
that Newton polygons can be used to define a stratification.
This stratification on $\A{g} \otimes {\FF}_p$  was 
much studied by Oort, who 
determined basic properties, see \cite{Oort-AM152}.

In order to define the Newton Polygon stratification, one takes for a
principally polarized abelian variety $X$ defined over 
$\overline{\FF}_q$ 
the Newton polygon of the action of geometric Frobenius $F_q$ on the 
cohomology group $H^1_{\rm et}(X \otimes \overline{\FF}_q, {\QQ}_{\ell})$,
or more generally the Newton polygon of the $p$-divisible group of $X$.
The Newton polygon is a symmetric polygon starting at $(0,0)$ and ending
at $(2g,g)$, lying below the line with slope $1/2$ and 
with integral vertices (break points). 
Here symmetric means that if slope $s$ occurs, then $1-s$ occurs
with the same multiplicity.

For principally polarized abelian varieties all such symmetric
polygons appear and thus there are 
$g(g-1)/2 +2$ strata in the NP stratification on $\A{g}\otimes {\FF}_p$.
The codimensions of the strata are known by Oort \cite{Oort-AM152}. 
The NP stratification depends only on the isogeny class of the abelian 
variety. 

The most degenerate NP stratum
is the supersingular stratum corresponding to slope $1/2$.
For $g=1$ and $g=2$ the $p$-rank zero stratum $V_0$ coincides with the
supersingular stratum, but not for $g\geq 3$. 

The second stratification is due to Ekedahl and Oort, 
see \cite{Oort1999, Oort2001}.
For the Ekedahl-Oort stratification one looks at the 
isomorphism type of the group scheme $X[p]$ together 
with Frobenius $F$ and Verschiebung $V$. 
Alternatively, one can
look at the de Rham cohomology $H^1_{\rm dR}(X)$ and at the relative position
of the kernels $\ker(F)$ and $\ker(V)$, as done in \cite{vdG1999}. 
These are totally isotropic subspaces
in the space $H^1_{\rm dR}(X)$, which is a non-degenerate symplectic space by
the Weil pairing. This stratification possesses $2^g$ strata. 
I showed in 1994  that the 
strata $Y_{\mu}$ are indexed by Young diagrams, or equivalently, by
tuples $\mu=[\mu_1,\ldots,\mu_r]$ with $0 \leq r \leq g$ and $\mu_i> \mu_{i+1}$. 
In fact, they can be interpreted as the degeneration loci of maps between
vector bundles, see \cite{vdG1999}.

The largest open stratum is the locus of ordinary abelian varieties. The codimension of $Y_{\mu}$ is $\sum_i \mu_i$.
The stratification can be extended to good toroidal compactifications of Faltings-Chai type. 

The Ekedahl-Oort and the Newton Polygon stratification share a number of strata:
the $p$-rank strata. The E-O strata are in general not preserved by isogenies
and deal with more subtle properties. 
The dimensions of the E-O strata are known. 
There has been a lot of research on this stratification, 
with a finer structure provided by a foliation, 
for example see \cite{Oort2004}.

\end{section}
%%%%%%%%%%%%%%%%%%%%%%%%%%%%%%%%%%%%
\begin{section}{Cycle Classes}
The most degenerate E-O stratum is the
superspecial locus. An abelian variety is superspecial if it isomorphic
to a product of supersingular elliptic curves (as an unpolarized abelian variety).
This superspecial locus has dimension $0$. All its points are rational over
${\FF}_{p^2}$. The stacky interpretation of Deuring's formula allows a
generalization. It was given by Ekedahl (\cite{Ekedahl1987}).

\begin{theorem} The number of superspecial abelian varieties is given by
$$
\sum_X \frac{1}{\# {\rm Aut}(X)} =
(p-1)(p^2-1) \cdots (p^g+(-1)^g) \, p(g)\, ,
$$
where the sum is over the isomorphism classes of principally polarized
superspecial abelian varieties $X$ over $\overline{\FF}_p$ and $p(g)$ is the
constant
$$
p(g)=(-1)^{g(g+1)/2} 2^{-g} \, \zeta(-1)\zeta(-3)\cdots \zeta(1-2g) \, ,
$$ 
where $\zeta$ denotes the Riemann zeta function.
\end{theorem}
This constant $p(g)$ is a proportionality constant related to the volume of
${\rm Sp}(2g,{\ZZ})\backslash \mathfrak{H}_g$. Here is a little table. 

\vbox{
\bigskip\centerline{\def\quad{\hskip 0.6em\relax}
\def\quod{\hskip 0.5em\relax }
\vbox{\offinterlineskip
\hrule
\halign{&\vrule#&\strut\quod\hfil#\quad\cr
height2pt&\omit&&\omit&&\omit&&\omit&&\omit&\cr
%\noalign{\hrule}
& $g$ && $1$ && $2$ && $3$ && $4$ &\cr
\noalign{\hrule}
& $p(g)$ && $1/24$ && $1/5760$  && $1/2903040$ && $1/1393459200$& \cr
} \hrule}
}}

\bigskip

The formula of Ekedahl can be interpreted as
 a formula for the degree of the class of the
$0$-dimensional superspecial stratum. As such it 
allows a far reaching generalization: the cycle classes
of the E-O strata are known. The interpretation of these strata 
as degeneracy classes of maps between 
vector bundles allows their calculation. 
We refer to \cite{vdG1999, E-vdG2009}.

For example, for the $p$-rank stratum $V_f$ on ${\mathcal{A}}_g \otimes {\FF}_p$ 
the cycle class is
$$
[V_f]=(p-1)(p^2-1)\cdots (p^{g-f}-1) \, \lambda_{g-f}\, ,
$$
where $\lambda_i$ is the $i$th Chern class of the Hodge bundle ${\EE}$. 
This formula holds on $\A{g}\otimes {\FF}_p$ but can be extended to 
good toroidal 
compactifications $\tilde{\mathcal{A}}_g \otimes {\FF}_p$. 
Such formulas can be seen as a generalization of Deuring's formula. 
Indeed, for $g=1$ 
the locus $V_0$ is the locus of supersingular elliptic curves and
it has class $(p-1)\lambda_1$. 
To connect it with Deuring's formula, observe that the modular
form $\Delta$, given in (0) in Section \ref{section-counting}, 
represents a section of ${\det}({\EE})^{12}$ and its 
zero divisor, which represents $12\, \lambda_1$,
is the cusp of $\tilde{\mathcal{A}}_1$, which 
represents a physical point with
multiplicity $1/2$ (since the degenerate elliptic curve it represents has
$\# {\rm Aut}=2$). Thus we find $\deg (V_0)=(p-1)/24$.

\bigskip

An intuitive explanation of the formula for the 
cycle class of
the $p$-rank locus $V_f$ may be given as follows.
If an abelian variety $X$ of dimension $g$ 
has $p$-rank $g$, then its $p$-kernel $X[p]$ 
contains the infinitesimal
group scheme $\mu_p^g$. Choosing a basis gives us $g$ tangent vectors
at the origin of $X$. The (open) $p$-rank locus $V_f-V_{f-1}$ is the locus where
$f$ sections survive. Chern classes measure the independence 
of sections. So it is natural that the $(g-f)$th Chern class of the Hodge bundle appears.

\bigskip

For the NP stratification the cycle classes of the strata are in
general not known for strata that do not occur as E-O strata. But 
the 
class of the supersingular locus $S_3$ for $g=3$, which is not
an E-O stratum, is known; its class is
$$
[S_3]=(p-1)^2(p^3-1)(p^4-1) \lambda_1 \lambda_3\, ,
$$ 
see \cite[Thm 11.3]{vdG1999}. We intend to come back to the
calculation of such cycle classes in the near future.

\end{section}
%%%%%%%%%%%%%%%%%%%5
\begin{section}{Strata on $\M{g}\otimes {\FF}_p$}

Via the Torelli morphism these two stratifications on $\A{g}\otimes {\FF}_p$ 
induce stratifications on $\M{g}\otimes {\FF}_p$. But here questions abound.
First: what are the dimensions of the strata? In particular: which strata are non-empty?

A first example is $g=2$. Here $\M{2}$ is an open subset of $\A{2}$. The cycle classes 
give information. To avoid the stacky aspect we may look at the moduli space
$\A2[n]$ of level $n$. This space is for $n\geq 3$ a variety and we can consider
for $p$ not dividing $n$ the cover
$$
\A{2}[n]\otimes {\FF}_p \langepijl{} \A{2}\otimes {\FF}_p \, .
$$ 
It is a Galois cover of
degree $r(n)=\# {\rm Sp}(4,{\ZZ}/n{\ZZ})$. 

The supersingular locus consists of a number of projective
lines, so-called Moret-Bailly lines. Indeed, every supersingular
principally polarized abelian surface can be obtained 
as a quotient $E^2/j(\alpha_p)$, where $E$ is a 
fixed supersingular elliptic curve over ${\FF}_p$, 
the product $E^2$ is provided with the
polarization Frobenius $F$, and where $j:\alpha_p\hookrightarrow
\alpha_p^2\cong E^2[F]$ 
is an embedding of the group scheme $\alpha_p$ into
the kernel $E^2[F]$ of $F$. These embeddings
are parametrized by a ${\PP}^1$ and it is easy to see that
the determinant of the Hodge bundle of the resulting family
over ${\PP}^1$ has degree $p-1$. We know the class of the 
supersingular locus $[V_0]=(p-1)(p^2-1)\lambda_2$ and since
the degree of $\lambda_1$ on each line is $p-1$  and 
$\deg(\lambda_1\lambda_2)= 1/5760$, we find
$(p^2-1) r(n)/5760$ projective lines.

We also know that the degree of the superspecial locus 
is $(p-1)(p^2+1)r(n)/5760$. This fits, since on each line we have
$p^2+1$ points that are rational over ${\FF}_{p^2}$ and so we
see that through
each superspecial point $p+1$ lines pass. 

So the supersingular locus
consists of $(p^2-1) r(n)/5760$ projective lines meeting in $(p-1)(p^2+1)r(n)/5760$
superspecial points. But only
$$
(p-1)(p-2)(p-3) \frac{r(n)}{5760}
$$
of these lie in $\M{2}[n]\otimes {\FF}_p$, as one sees by subtracting 
those that lie in the complement  $\A{1,1}[n]\otimes {\FF}_p$. 
This simple example 
shows the use of the cycle classes and agrees with results in 
\cite{K-O}. 

In particular, we see that there are no such superspecial points 
on $\M{2}\otimes {\FF}_p$ for $p=2$ and $p=3$.  
Of course, this is related to the well-known fact that there are no maximal curves
of genus $2$ over ${\FF}_4$ and ${\FF}_9$. 

The formula just given  illustrates 
that strata can be empty on 
$\M{g}\otimes {\FF}_p$. But the $p$-rank strata on $\M{g}\otimes {\FF}_p$ have
the expected dimension as Faber and I showed.

\begin{theorem} {\rm (\cite{F-vdG2})}
For every $p$ the locus in $\M{g}\otimes {\FF}_p$ of curves with $p$-rank $\leq f$
is pure of dimension $g-f$.
\end{theorem}
The same idea can be applied to other cases, see \cite{A-P2}.
Recently, there has been a lot of activity trying to investigate the 
dimensions of NP and E-O strata on $\M{g}\otimes {\FF}_p$ for low values
of $g$ by constructions of explicit curves or families of curves over
finite fields, see for example papers by Rachel Pries and others, \cite{A-P1, Pries1, LMPT1,LMPT2}. 
But our knowledge is still very limited.
\end{section}
%%%%%%%%%%%%%%%%%%%%%%%%%5
\begin{section}{Supersingular curves}
Starting at the other end of the NP stratification  on $\M{g}\otimes {\FF}_p$
we may ask: Is there a supersingular curve of genus $g$ in characteristic $p$?
Maximal curves over a finite field are supersingular, so if available 
provide an answer. For example, Ibukiyama showed the existence for $g=3$
and $p>2$.

\begin{theorem} {\rm (\cite{Ibukiyama})}
For odd $p$ there exists a curve of genus $3$ over ${\FF}_p$ whose Jacobian is
superspecial over ${\FF}_{p^2}$.
\end{theorem}  
It follows from Ekedahl's formula 
that there is no superspecial curve of genus $3$ for $p=2$.

For $g=4$ we have a positive answer by Kudo, Harashita and Senda.

\begin{theorem} {\rm (\cite{K-H-S})}
There exists a supersingular curve of genus $4$ in every characteristic.
\end{theorem}

But one can also fix $p$ and vary $g$. Coding theory led 
van der Vlugt and me to consider curves 
related to Reed-Muller codes in \cite{vdG-vdV1992a}. 
For $q$ a power of $2$ these curves over
${\FF}_q$ are given by an equation
$$
y^2+y=xR(x)
$$
with $R=\sum_{i=0}^h a_i x^{2^i}\in {\FF}_q[x]$  
a so-called $2$-linearized polynomial. 
These curves are supersingular of genus $2^{h-1}$ for $R$ of degree $2^h$.
They possess large extra-special groups of automorphisms. 
Using such curves as building blocks and using fibre products 
one can construct supersingular 
curves of arbitrary genus in characteristic $2$, even over ${\FF}_2$.

\begin{theorem} {\rm (\cite{vdG-vdV1992a})}
For every $g$ there exists a supersingular curve of genus $g$ 
defined over ${\FF}_2$.
\end{theorem}
To give one example, for $g=2021$ we have the supersingular curve over ${\FF}_2$
$$
y^{256}+y^{64}+y^4+y= x^{68}+x^{20}+x^{17}+x^{12}+x^{10} \, .
$$
The same method can be applied to the case $p>2$ with as building blocks
Artin-Schreier curves $y^p-y=xR(x)$ with $R=\sum_{i=0}^h a_ix^{p^i}$
a $p$-linearized polynomial of degree $p^h$ of genus $p^h(p-1)/2$, as studied in
\cite[Section 13]{vdG-vdV1992a}, but here these do not cover all genera. 

\begin{proposition}
For odd $p$ there exists a supersingular curve over ${\FF}_p$ 
for every genus g with only $0$ and $(p-1)/2$ in its $p$-adic expansion.
\end{proposition}
As an example, for $g=999$ one finds over ${\FF}_3$ the supersingular curve
$$
y^{27}+y^9+y^3+y=x^{246}+x^{84}+x^{82} \, .
$$
Unfortunately, this proposition covers only a very thin set of genera 
for $p>2$. But one can do variations. 
One can take quotients of these curves, or of other supersingular
curves.
For example, if $d$ is a divisor of $p^h+1$
then the curve $y^{p^m}-y=x^d$ for $m\geq 1$ and $h\geq 0$ 
is supersingular because it can be obtained
as a quotient of a curve $y^{p^m}-y=x^{p^h+1}$.

\smallskip

For $3 \leq p \leq 23$ Riccardo Re constructed a supersingular curve in 
characteristic $p$ for almost all $g\leq 100$ with very few undecided exceptions,
see \cite{Re1}. He used a variety of methods, for examples taking quotients
of Fermat type curves. For a related reference see \cite{Bouw}.

\end{section}
%%%%%%%%%%%%%%%%%%%%%%%
\begin{section}{Bounds on the $a$-number}
Besides the $p$-rank of an abelian variety there is another invariant
closely related to the E-O stratification, the $a$-number, 
introduced by Oort.
It measures the dimension of the intersection of
$\ker(F)$ and $\ker(V)$ in $X[p]$ (or in $H^1_{\rm dR}(X)$) 
and may be defined by
$$
a(X)=\dim_k \left(\ker F^{\ast} : H^1(X,O_X) \to H^1(X,O_X)\right)\, .
$$
We have $0 \leq a(X) \leq g$. For a curve $C$ we put $a(C)=a({\rm Jac}(C))$ and
have
$$
a(C)=g - {\rm rank}(V)\, ,
$$
with $V: H^0(C,\Omega_C^1)\to H^0(C,\Omega_C^1)$ the Cartier operator.
The loci $T_a$ of abelian varieties with $a$-number $\geq a$ are (closed) 
strata of the E-O stratification on $\A{g}\otimes {\FF}_p$ 
and have codimension $a(a+1)/2$. 

The most special case of the E-O stratification is the case of superspecial 
abelian varieties. 
Oort showed in 1975 \cite{Oort1975}: an abelian variety $X$ is superspecial 
if and only if $a(X)=g$.

For superspecial Jacobians of curves  
there is a theorem of Ekedahl 
that limits the genus
of a curve with a superspecial Jacobian.

\begin{theorem}\label{Ekedahl-Thm} {\rm (\cite{Ekedahl1987})}
If $C$ is a superspecial curve of genus $g$ in characteristic $p$ then $g \leq p(p-1)/2$\, .
\end{theorem}
We thus see that for large $g$ the superspecial stratum is empty 
on $\M{g}\otimes {\FF}_p$.
We can interpret Ekedahl's result Theorem \ref{Ekedahl-Thm} 
as: $a=g$ implies $g\leq p(p-1)/2$.
Zijian Zhou, improving work of Riccardo Re (\cite{Re2}),
showed an upper bound on the genus for the
case that $a=g-1$.

\begin{theorem} {\rm (\cite{Zhou})}
If $C$ is a curve in characteristic $p$ with $a(C)=g-1$, then $g \leq p+ p(p-1)/2$.
\end{theorem}

One may ask for a generalization.
Here is my conjecture for the optimal result.

\begin{conjecture}
For a curve $C$ of genus $g$ in characteristic $p>0$ we have
$$
a(C) \leq \frac{p-1}{p} g + \frac{p-1}{2}\, ,
$$
equivalently, with $a(C)=g-r$ with $r$ the rank of the Cartier operator,
we have
$$
g \leq p\, r +\frac{p(p-1)}{2} \, .
$$
\end{conjecture}
Note that this is in accordance with the results of Ekedahl and Zhou.
Moreover, for $p=2$ this gives $a\leq (g+1)/2$, 
a result of St\"ohr-Voloch, \cite{S-V}. They
show that this inequality is strict if 
$g\geq 3$ and $C$ non-hyperelliptic.

The conjecture predicts empty strata of the E-O stratification on $\M{g}\otimes {\FF}_p$.
\end{section}
%%%%%%%%%%%%%%%%%%%%%%%%%%%%%%%%%%%%
\begin{section}{Counting points on strata}
An interesting question is: which strata on $\M{g}\otimes {\FF}_p$ 
have dimension $0$? For these strata one should determine
their number of rational points. 
It could provide interesting curves over finite fields.

\bigskip

One can try to count number of points over finite fields on strata. 
Here the foliations introduced by Oort should play a role. This could
lead to congruences for modular forms.
In general not much is known. 

A simple example are the E-O strata on 
$\M{2}\otimes {\FF}_2$. 
The number of points $\# S({\FF}_q)$ on a (closed) stratum $S$ for $q=2^m$
were given in \cite{vdG-vdV1992b} and are shown
in the following table. We index the strata 
by the $2$-rank $f$ 
or by the $a$-number.

\vskip 0.5 cm

\begin{center}
\begin{tabular}{|l|r|r|r|r|r|}
\hline
stratum & $f\leq 2$ & $f\leq 1$ & $f=0 $ & $a=2$  \\
\hline
% & & & &   \\
$\# {S}({\FF}_q) $ & $q^3$ & $q^2$  & $q$ & $0$ \\
\hline
\end{tabular}
\end{center}

\vskip 0.5 cm 

One may interpret some of the results of \cite{K} on the supersingular locus for $g=3$ as a result
in this vein. 

Nart and Ritzenthaler (see \cite{N-R}) gave the cardinalities of 
${\FF}_q$-rational points for $q=2^m$ 
for the NP strata on $\M{3}\otimes {\FF}_2$. 
We give a table with their results, where
the strata $S$ are indicated by the $p$-rank $f$ or the slope $s=1/2$.
Also the cardinalities of the intersection with the hyperelliptic locus
$\mathcal{H}_3$ are given.

\vskip 0.5 cm

\begin{center}
\begin{tabular}{|l|r|r|r|r|r|r|}
\hline
stratum & $f\leq 3$ & $f\leq 2$ & $f\leq 1$ & $f=0$  & $s=1/2$ \\
\hline
$\# {S}({\FF}_q) $ & $q^6+q^5+1$ & $q^5+q^4$  & $q^4+2q^3-q^2$ & $q^3+q^2$ & $q^2$ \\
$\#({S}\cap \mathcal{H}_3)({\FF}_q) $ & $q^5$ & $q^4$  & $2q^3-q^2$ & $q^2$ & $0$ \\
\hline
\end{tabular}
\end{center}

\bigskip

The author hopes that this text will entice some readers to engage in
the many possibilities for fruitful and pleasant 
experimentation in this fascinating corner of mathematics. 
\end{section}
%%%%%%%%%%%%%%%%%%%%%%%%%%%%%%%%%%%%%%%%


\begin{thebibliography}{99}
%%%%%%%%%%%%%%%%%%%%%%%%%%%%%%%%%%%%%%%%
\bibitem{A-P1} J.\ Achter, R.\ Pries:
{\sl Monodromy of the $p$-rank strata of the moduli space 
of curves.} Int.\ Math.\ Res.\ Not.\ IMRN 2008, no. 15.

\bibitem{A-P2} J.\ Achter, R.\ Pries: 
{\sl The $p$-rank strata of the moduli space of 
hyperelliptic curves.}  Adv.\ Math.\ \textbf{227} (2011), 
1846–-1872.

\bibitem{Behrend} K.\ Behrend: 
{\sl The Lefschetz trace formula for algebraic stacks.}
Inv.\ Math.\ \textbf{112} (1993), 127--149.

\bibitem{JB1} J. Bergstr\"om: {\sl Cohomology of moduli spaces of curves
of genus three via point counts.}
J. Reine Angew. Math. \textbf{622} (2008), 155--187.

\bibitem{JB2} J. Bergstr\"om: {\sl Equivariant counts of points on the moduli
spaces of hyperelliptic curves.} Documenta Math.\ \textbf{14} (2009), 259--296.


\bibitem{BFvdG} J.\ Bergstr\"om, C.\ Faber, G.\ van der Geer:
{\sl Siegel modular forms of genus 2 and level 2:
cohomological computations and conjectures.}
Int.\ Math.\ Res.\ Not.\ IMRN (2008), Art.\ ID rnn 100, 20 pp.

\bibitem{BFvdG2} J.\ Bergstr\"om, C.\ Faber, G.\ van der Geer:
{\sl Siegel modular forms of degree three and
the cohomology of local systems.}
Selecta Math. (N.S.) \textbf{20} (2014), no. 1, 83--124.

\bibitem{BFvdG3} J.\ Bergstr\"om, C.\ Faber, G.\ van der Geer:
Siegel Modular Forms of Degree Two and Three. 2017. Website.
{\tt http://smf.compositio.nl}

\bibitem{B-vdG} J.\ Bergstr\"om, G.\ van der Geer:
{\sl Picard modular forms and the cohomology of local systems on a
Picard modular surface.}  
Comment.\ Math.\ Helv.\ \textbf{97} (2022), 305--381.

\bibitem{B-T} J.\ Bergstr\"om, O.\ Tommasi:
{\sl The rational cohomology of $\barM{4}$.}
Math.\ Annalen \textbf{338} (2007), 207--239.

\bibitem{Birch}
B.J.\ Birch: {\sl How the number of points of an elliptic
curve over a fixed prime field varies.}
Journal  London Math.\ Soc., \textbf{43} (1968), p.\ 57--60.

\bibitem{vdB-E2005} T.\ van den Bogaart, B.\ Edixhoven:
{\sl Algebraic stacks whose number of points over finite fields is a polynomial.}
In: Number fields and functions fields--two parallel worlds. 
Progress in Math.\ 239. Birkh\"auser, 2005, p.\ 39--49.

\bibitem{Bouw} I.\ Bouw, W.\ Ho, B.\ Malmskog, R.\ Scheidler, P.\ Srinivasan, C.\ Vincent:
{\sl Zeta functions of a class of Artin-Schreier curves with many automorphisms.}
Directions in number theory, 87--124. Assoc.\ Women Math.\ Ser., 3, Springer 

\bibitem{BGHZ} J.\ Bruinier, G.\ van der Geer, G.\ Harder, D.\ Zagier:
{\sl The 1-2-3 of modular forms.} Universitext. Springer Verlag 2007.

\bibitem{Chenevier-Lannes} G.\ Chenevier, J.\ Lannes: {\sl Automorphic forms and even unimodular lattices. Kneser neighbors of Niemeier lattices.} Translated from the French by Reinie Ern\'e. Ergebnisse der Mathematik und ihrer Grenzgebiete (3), 69. Springer, Cham, 2019. xxi+417 pp.

\bibitem{Chenevier-Renard}  G.\ Chenevier, D.\ Renard: 
{\sl Level one algebraic cusp forms of classical groups of small rank.} Mem.\ Amer.\ Math.\ Soc.\ \textbf{237} (2015), no. 1121, v+122 pp.


\bibitem{CFvdG} F.\ Cl\'ery, C.\ Faber, G.\ van der Geer:
{\sl Concomitants of ternary quartics and 
vector-valued Siegel and Teichmüller modular forms 
of genus three.} Selecta Math.\ (N.S.) \textbf{26} (2020), 
no. 4, Paper No.\ 55, 39 pp.

\bibitem{Deligne1968} P.\ Deligne: 
{\sl Formes modulaires et repr\'esentations $\ell$-adiques.}
S\'eminaire N.\ Bourbaki, 1968--1969, exp.\ 355, p.\ 139--172.

\bibitem{Deligne1974}  P.\ Deligne:
{\sl Th\'eorie de Hodge III.} Publ.\ Math.\ IHES \textbf{44} (1974),
5--77.

\bibitem{Deligne} P.\ Deligne:
{\sl Cohomologie \'etale.} SGA $4\frac{1}{2}$.
Lecture Notes in Math.\ \textbf{569}, Springer Verlag, 1977.

\bibitem{Deligne-Mumford} P.\ Deligne, D.\ Mumford:
{\sl The irreducibility of the spaceof curves of given genus.}
Publ.\ Math.\ de I.H.E.S. \textbf{36}, (1969), 75--109.

\bibitem{Diaconu2020} 
A.\ Diaconu:
{\sl Equivariant Euler characteristics of $\barM{g,n}$.}
Algebraic Geometry \textbf{7} (2020), 523–-543.

\bibitem{Ekedahl1987}  T.\ Ekedahl: 
{\sl On supersingular curves and abelian varieties.}
Math.\ Scand.\ \textbf{60} (1987), 151--178.

\bibitem{E-vdG2009} T.\ Ekedahl, G.\ van der Geer:
{\sl Cycle classes of the E-O stratification and 
the moduli of abelian varieties.}
In: Algebra, Arithmetic and Geometry. (Y.\ Tschinkel, Y.\ Zarhin eds.)
Progress in Math.\ 269, 567--636. Birkh\"auser 2009.

\bibitem{F-vdG1} C.\ Faber, G.\ van der Geer:
{\sl Sur la cohomologie des syst\`emes locaux sur les 
espaces des modules des courbes de genre $2$ et des 
surfaces ab\'eliennes. I, II.} 
C.R.\ Acad.\ Sci.\ Paris, S\'er. I, \textbf{338} (2004), 381–-384,
467–-470.

\bibitem{F-vdG2} C.\ Faber, G.\ van der Geer:
{\sl Complete subvarieties of $\M{g}$ and the Prym map.}
Journal f\"ur die reine und angewandte Mathematik \textbf{573} (2004), 117–-137.
 
\bibitem{F-G1} X.\ Faber, J.\ Grantham:
{\sl Binary curves of small fixed genus and gonality with many fixed points}
 Journal of Algebra \textbf{597} (2022), 24--46.

\bibitem{F-G2} X.\ Faber, J.\ Grantham:
{\sl Ternary and quaternary curves of small fixed genus and gonality with many rational points.} {\tt arXiv:2010.07992}

\bibitem{F-G-H}  X.\ Faber, J.\ Grantham, E.\ Howe:
{\sl On the maximum gonality of a curve over a finite field.}
{\tt arXiv:2207.14307}
\bibitem{F-C} G.\ Faltings, C.-L.\ Chai:
{\sl Degeneration of abelian varieties.}
Ergebnisse der Mathematik und ihrer Grenzgebiete (3), 22.
Springer-Verlag, Berlin, 1990.

\bibitem{F-H} W.\ Fulton, J.\ Harris:
{\sl Representation theory. A first course.} 
Graduate Texts in Mathematics, 129. 
Readings in Mathematics. 
Springer-Verlag, New York, 1991. 

\bibitem{vdG1999} G.\ van der Geer: 
{\sl Cycles on the moduli space of Abelian varieties.}
In: Moduli of Curves and Abelian varieties. (C.\ Faber, E.\ Looijenga eds) 
The Dutch Intercity Seminar. Aspects of Math.\ 33 (1999), 65--89.

\bibitem{vdG2011} G.\ van der Geer: 
{\sl Rank one Eisenstein cohomology of local systems on
the moduli space of abelian varieties.}
Science in China, Math.\ \textbf{54}, 1621--1634 (2011).

\bibitem{vdG-ECM} G.\ van der Geer:
{\sl  Curves over finite fields and codes.}
In: European Congress of Mathematics, Vol.\ II (Barcelona, 2000),
Progr.\ Math.\ 202,  225-–238. Birkh\"auser, Basel, 2001.
\bibitem{Getzler1998} E.\ Getzler:
{\sl The semi-classical approximation for modular operads.}
Comm.\ Math.\ Phys.\ \textbf{194} (1998), 481–-492.

\bibitem{vdG-vdV1992a} G.\ van der Geer, M.\ van der Vlugt:
{\sl Reed-Muller codes and supersingular curves I.}
Compositio Mathematica \textbf{84} (1992), 333--367.

\bibitem{vdG-vdV1992b} G.\ van der Geer, M.\ van der Vlugt:
{\sl Supersingular curves of genus $2$ over finite fields 
of characteristic~$2$.} Mathematische Nachrichten \textbf{159} (1992), 73–-81.

\bibitem{vdG-vdV1995} G.\ van der Geer, M.\ van der Vlugt:
{\sl On the existence of supersingular curves of given genus.}
Journal f\"ur die reine und angewandte Mathematik \textbf{458} 
(1995), 53--61.

\bibitem{Getzler1998b} E.\ Getzler:
{\sl Topological recursion relations in genus $2$.}
In: Integrable systems and algebraic geometry. (Kobe/Kyoto 1997),
p.\ 73--106.
World Sci.\ Publishing, River Edge, NJ, 1998, 

\bibitem{Getzler1999} E.\ Getzler: 
{\sl Resolving mixed Hodge modules on configuration spaces.}
Duke Mathematical Journal \textbf{96} (1999), 175--203.

\bibitem{Getzler2002} E.\ Getzler:
{\sl Euler characteristics of local systems on $\M{2}$.}
Compositio Mathematica \textbf{132} (2002), 121--135.

\bibitem{G-K} E.\ Getzler, M.\ Kapranov:
{\sl Modular operads.} Compositio Mathematica  \textbf{110} (1998), 65--126.

\bibitem{Harder}
G.\ Harder:
{\sl The Eisenstein motive for the cohomology of 
${\rm GSp}_2({\ZZ})$,} pp.\ 143--164, in: Geometry and arithmetic 
(Island of Schiermonnikoog, 2010), edited by C.\ Faber 
et al., European Mathematical Society, Z\"urich, 2012.

\bibitem{Homma} M.\ Homma:
{\sl Automorphisms of prime order of curves.}
Manuscripta Mathematica \textbf{33} (1980), 99--109.

\bibitem{Ibukiyama} T.\ Ibukiyama:
{\sl On rational points of curves of genus $3$ over finite fields.}
Tohoku Mathematical Journal \textbf{45} (1993), 311--329

\bibitem{Ichikawa} T.\ Ichikawa:
{\sl Teichm\"uller modular forms of degree $3$.}
American  Journal of  Mathematics \textbf{117}
(1995), 1057--1061.

\bibitem{Igusa1962} J.-I.\ Igusa: 
{\sl On Siegel modular forms of genus two.}
 Amererican Journal of  Mathematics \textbf{84} (1962), 175--200.

\bibitem{Igusa1967} J.-I.\ Igusa:
{\sl Modular forms and projective invariants.}
American Journal of Mathematics \textbf{89} (1967), 817--855.

\bibitem{K} V.\ Karemaker, F.\ Yokubo, C-F.\ Yu:
{\sl Mass formula and Oort's conjecture for supersingular abelian threefolds.} 
Advances in Mathematics \textbf{386}, (2021). 

\bibitem{K-O} T.\ Katsura, F.\ Oort: 
{\sl Families of supersingular abelian surfaces.} 
Compositio Mathematica \textbf{ 62} (1987), 107--167.

\bibitem{K-L} M.\ Kisin, G.I.\ Lehrer:
{\sl Equivariant Poincar\'e polynomials and counting points over finite fields.}
Journal of Algebra \textbf{247}, (2002),  435–-451.

\bibitem{Koblitz} N.\ Koblitz:
{\sl $p$-adic variation of the zeta-function over families of
varieties defined over finite fields.}
Compositio Mathematica, \textbf{31} (1975),  119--218.

\bibitem{K-H-S} M.\ Kudo, S.\ Harashita, H.\ Senda:
{\sl The existence of supersingular curves of genus $4$ in
arbitrary characteristic.}
Research in Number Theory (2020).
  
\bibitem{LMPT1} W.\ Li, E.\ Mantovan, R.\ Pries, Y.\ Tang:
{\sl Newton polygon stratification of the Torelli locus in PEL-type 
Shimura varieties.} {\tt arXiv:1811.00604}

\bibitem{LMPT2} W.\ Li, E.\ Mantovan, R.\ Pries, Y.\ Tang:
{\sl Newton polygons arising for special families of cyclic 
covers of the projective line.} Res.\ Number Theory \textbf{5} (1), 12 (2019).

\bibitem{Moonen} B.\ Moonen: 
{\sl Special subvarieties arising from families of cyclic covers of the projective line.} 
Documenta Mathematica \textbf{15}, 793--819 (2010).

\bibitem{Mumford} D.\ Mumford:
{\sl Geometric invariant theory.} 
Ergebnisse der Math.\ \textbf{34} Springer Verlag, 1965.

\bibitem{N-R} E.\ Nart, C.\ Ritzenthaler: 
{\sl Jacobians in isogeny classes of supersingular 
abelian varieties in characteristic $2$.} 
Finite Fields Appl.\ \textbf{14}, No. 3, 676--702 (2008).

\bibitem{Oort1975} F.\ Oort: 
{\sl Which abelian surfaces are products of
elliptic curves?} Mathematische Annalen \textbf{214} (1975), 35--47.

\bibitem{Oort1999} F.\ Oort: 
{\sl A stratification of a moduli space of polarized abelian varieties
in positive characteristic.}
In: Moduli of curves and abelian varieties. The Dutch Intercity Seminar.
(C.\ Faber, E.\ Looijenga eds)
Aspects of Mathematics \textbf{33} (1999), 
 
\bibitem{Oort-AM152} F.\ Oort: 
{\sl Newton polgons and formal groups.}
Annals of Mathematics \textbf{152} (2000), 183--206.

\bibitem{Oort2001} F.\ Oort:
{\sl A stratification of a moduli space of abelian varieties.}
 Progress in Mathematics \textbf{195}, (2001), 345-–416, Birkh\"auser,
Basel, (2001)

\bibitem{Oort2004} F.\ Oort:
{\sl Foliations in moduli spaces of abelian varieties.}
Journal of the American Mathematical Society \textbf{17}, No. 2, 267--296 (2004).

\bibitem{Petersen} D.\ Petersen:
{\sl Cohomology of local systems on the moduli of principally polarized abelian surfaces.}
Pacific Journal of  Mathematics \textbf{ 275} (2015), no. 1, 39-–61. 

\bibitem{Pries1} R.\ Pries: 
{\sl Current results on Newton polygons of curves.} In: 
Open problems in arithmetic algebraic geometry. Adv.\ Lect.\
Math.\ \textbf{46} (2019), 179–-207 (2019).
Int.\ Press, Somerville, MA. 

\bibitem{Re1} R.\ Re: 
{\sl Invariants of curves and Jacobians in positive characteristic.}
Ph.D. thesis. University of Amsterdam 2004.

\bibitem{Re2} R.\ Re: 
{\sl the rank of the Cartier operator and linear systems on curves.}
Journal of Algebra \textbf{236} (2001), 80--92.

\bibitem{Ritzenthaler} C.\ Ritzenthaler: {\sl Postface to Chapter 4.}
In: J.P. Serre: Rational points on curves over finite fields.
Documents Math\'ematiques 18. SMF 2020, pages 84--90.

\bibitem{Rohde} J.\ Rohde:
{\sl 
Cyclic coverings, Calabi-Yau manifolds and complex multiplication.}
Lecture Notes in Math.\ \textbf{1975}, (2000). Springer Verlag.

\bibitem{Roesner} M.\ R\"osner:
{\sl Parahoric restriction for ${\rm GSP}(4)$ and the inner 
cohomology of Siegel modular threefolds.} Inaugural dissertation, 
Ruprecht-Karls-Universität Heidelberg, 2016.

\bibitem{Satoh} T.\ Satoh: 
{\sl On certain vector valued Siegel modular forms of degree 2.}
Mathematische Annalen \ \textbf{274} (1986), 335–-352.

\bibitem{Scholl} A.J.\ Scholl: {\sl
Motives for modular forms.}  Inventiones Math.\
\textbf{100}  (1990),  no. 2, 419--430.

\bibitem{Serre} J.P.\ Serre: 
{\sl Rational points on curves over finite fields.}
Documents Math\'ematiques 18. SMF 2020.

\bibitem{Shimura1959} G.\ Shimura:
{\sl Sur les int\'egrales attach\'ees aux formes automorphes.}
J.\ Math.\ Soc.\ Japan \textbf{11}, (1959), 291--311.

\bibitem{Shimura1964} G.\ Shimura: {\sl On purely transcendental fields of
automorphic functions of several variables.} {Osaka Math.\ Journal \bf 1}
(1964), 1--14.

\bibitem{S-V} K.-O.\ St\"ohr, J.F.\ Voloch:
{\sl A formula for the Cartier operator on plane algebraic curves.} 
Journal f\"ur die reine und angewandte Mathematik \textbf{377} (1987), 49–-64.
 
\bibitem{Tommasi} O.\ Tommasi:
{\sl Rational cohomology of the moduli space of genus $4$ curves.}
Compositio Mathematica \textbf{141} (2005), 359--384.

\bibitem{Vermeulen} F.\ Vermeulen: 
{\sl Curves of fixed gonality with many rational points.}
{\tt arXiv:2102.00900}
 

\bibitem{Weil} A.\ Weil:
{\sl The field of definition of a variety.}
American Journal of  Mathematics \textbf{78} (1956), 509-–524.

\bibitem{Weissauer} R.\ Weissauer:
{\sl The trace of Hecke operators on the space of classical holomorphic Siegel modular forms of genus two.}
{\tt arXiv:0909.1744}.

\bibitem{Xarles} X.\ Xarles:
{\sl A census of all genus $4$ curves over the field with $2$ elements.}
arXiv:2007.07822v1.

\bibitem{Zhou} Z.\ Zhou: 
{\sl The $a$-number and the Ekedahl-Oort types of Jacobians of curves.}
Ph.D. thesis University of Amsterdam 2019.
 

\end{thebibliography}
\end{document}